\documentclass[3p,fleqn,sort]{elsarticle}
\journal{Journal of Computational Physics}
\usepackage[labelfont=bf, singlelinecheck=false,
            justification=justified]{caption}
\usepackage{xcolor}
\usepackage{amsmath}
\usepackage{amssymb}
\usepackage{mathrsfs}
\usepackage{stmaryrd} 
\usepackage{booktabs}
\usepackage{microtype}
\usepackage{nicefrac}
\usepackage{siunitx}
\usepackage{mathtools} 
\usepackage{multirow}
\usepackage{enumerate}
\usepackage{subfig}
\usepackage{fixltx2e}
\usepackage{setspace}
\usepackage{placeins} 
\usepackage{etoolbox}
\usepackage{fancybox}
\robustify\bfseries

\usepackage[makeroom]{cancel}

\renewcommand{\vec}[1]{\ensuremath{\boldsymbol #1}}
\newcommand{\mat}[1]{\ensuremath{\mathbf{#1}}}
\newcommand{\avg}[1]{\left\{\!\left\{#1\right\}\!\right\}}

\newcommand{\half}{\frac{1}{2}}

\newcommand{\mmat}{\mathbf{M}}
\newcommand{\dmat}{\mathbf{D}}
\newcommand{\qmat}{\mathbf{Q}}

\newcommand{\comment}[1]{{\textcolor{black}{#1}}}

\begin{document}

\begin{frontmatter}

\title{A comparative study on polynomial dealiasing and split form discontinuous Galerkin schemes for under-resolved turbulence computations}

\author[mathematik]{Andrew R.~Winters\corref{correspondingauthor}}
\cortext[correspondingauthor]{Corresponding author}
\ead{awinters@math.uni-koeln.de}
\author[aeronautics,brasil]{Rodrigo C.~Moura}
\author[sciences]{Gianmarco Mengaldo}
\author[mathematik]{Gregor J.~Gassner}
\author[physik]{Stefanie Walch}
\author[aeronautics]{Joaquim Peiro}
\author[aeronautics]{Spencer J.~Sherwin}
\address[mathematik]{Mathematisches Institut, Universit\"at zu K\"oln, Weyertal 86-90, 50931 K\"oln, Deutschland}
\address[physik]{I.\,Physikalisches Institut, Universit\"at zu K\"oln, Z\"ulpicher Stra\ss{}e~77, 50937 K\"oln, Deutschland}
\address[aeronautics]{Department of Aeronautics, Imperial College London, South Kensington Campus, London SW7 2AZ, United Kingdom}
\address[sciences]{Division of Engineering and Applied Science, California Institute of Technology, Pasadena, California 91125, U.S.A.}
\address[brasil]{Depto.\ de Aerodin\^amica, Instituto Tecnol\'ogico de Aeron\'autica (ITA), Pra\c{c}a Mal.\ Eduardo Gomes, 50, S.\ J.\ Campos, Brasil}

\numberwithin{equation}{section}

\begin{abstract}
This work focuses on the accuracy and stability of high-order nodal discontinuous Galerkin (DG) methods for under-resolved turbulence computations. In particular we consider the inviscid Taylor-Green vortex (TGV) flow to analyse the implicit large eddy simulation (iLES) capabilities of DG methods at very high Reynolds numbers. The governing equations are discretised in two ways in order to suppress aliasing errors introduced into the discrete variational forms due to the under-integration of non-linear terms. The first, more straightforward way relies on consistent/over-integration, where quadrature accuracy is improved by using a larger number of integration points, consistent with the degree of the non-linearities. The second strategy, originally applied in the high-order finite difference community, relies on a split (or skew-symmetric) form of the governing equations. Different split forms are available depending on how the variables in the non-linear terms are grouped. The desired split form is then built by averaging conservative and non-conservative forms of the governing equations, although conservativity of the DG scheme is fully preserved. A preliminary analysis based on Burgers' turbulence in one spatial dimension is conducted and shows the potential of split forms in keeping the energy of higher-order polynomial modes close to the expected levels. This indicates that the favourable dealiasing properties observed from split-form approaches in more classical schemes seem to hold for DG. The remainder of the study considers a comprehensive set of (under-resolved) computations of the inviscid TGV flow and compares the accuracy and robustness of consistent/over-integration and split form discretisations based on the local Lax-Friedrichs and Roe-type Riemann solvers. Recent works showed that relevant split forms can stabilize higher-order inviscid TGV test cases otherwise unstable even with consistent integration. Here we show that stable high-order cases achievable with both strategies have comparable accuracy, further supporting the good dealiasing properties of split form DG. The higher-order cases achieved only with split form schemes also displayed all the main features expected from consistent/over-integration. Among test cases with the same number of degrees of freedom, best solution quality is obtained with Roe-type fluxes at moderately high orders (around sixth order). Solutions obtained with very high polynomial orders displayed spurious features attributed to a sharper dissipation in wavenumber space. Accuracy differences between the two dealiasing strategies considered were, however, observed for the low-order cases, which also yielded reduced solution quality compared to high-order results.
\end{abstract}

\begin{keyword}
Spectral element methods \sep Discontinuous Galerkin \sep Polynomial dealiasing \sep Split form schemes \sep Implicit large eddy simulation \sep Inviscid Taylor-Green vortex
\end{keyword}

\end{frontmatter}

\section{Introduction}
\label{intro}

In computational fluid dynamics (CFD), high-fidelity simulations have been recognised as a crucial area of investigation for understanding the physics underlying high Reynolds number turbulent flows, thereby providing the key tools for off-cycle aerodynamic design by means of physics-based simulations \cite{Pi2013}. However, current numerical schemes adopted in industry might not be suitable to provide an effective solution for these problems due to a highly dissipative and dispersive behaviour that damps and distorts turbulent structures, thus altering the energy transfer across scales and ultimately corrupting flow physics. An attractive class of numerical methods that is now widely spread in academia and is perhaps taking its first steps into industry is constituted by the so-called spectral element methods (SEM). Among these, one has the well-known continuous and discontinuous Galerkin (CG/DG) methods \cite{KarniadakisSherwin,hesthaven2007}, but also more recently developed formulations such as spectral volume \cite{wang2002spectral}, spectral difference \cite{liu2006spectral} and flux reconstruction \cite{hyunh} schemes. The present work focuses on the DG method, but its results may also be of special interest to those working with flux reconstruction, as the latter is able to recover certain variants of DG \cite{de2014connections,mengaldo2015discontinuous,mengaldo2016connections}.

Spectral element methods combine the advantageous properties of finite element/finite volume and spectral methods, namely geometric flexibility and reduced dispersion/dissipation errors (see e.g.\ \cite{mengaldo2017spatial,mengaldo2018spatial,fernandez2018non} on the latter topic). Nevertheless, the application of spectral element methods to challenging problems such as turbulent flows over complex geometries is still somewhat hindered by numerical instability issues. These are twofold: (a) the overall dissipation of high-order spectral element approaches might be insufficient to stabilize high Reynolds number turbulence computations, which are normally under-resolved due to computational cost constraints and therefore have significant energy at the smaller scales captured \cite{tcfd2012,moura2016ontheeddy}; (b) the most efficient spectral element methods that help alleviate the cost requirements commonly rely on under-integration of the non-linear terms, see e.g.\ \cite{Hindenlang2012}, which induces aliasing-driven instabilities that can also lead to numerical divergence \cite{Kirby2003}. If the discretisation order is not very high, dissipative effects might overcome polynomial aliasing errors and suppress instabilities, but accuracy is often compromised. At higher orders, however, aliasing and (lack of) dissipation issues require careful consideration and sometimes robustness may not be guaranteed even with dealiasing \cite{mengaldo2015,moura2016ontheeddy,kopriva2017aliasing}.

Issue (a) above can be addressed through the use of (dissipative) subgrid-scale models, in what would be effectively a large-eddy simulation (LES) methodology. Promising candidates within this category are based on the so-called variational multiscale (VMS) approach introduced by Hughes et al.\ \cite{hughes1998variational,hughes2000large}. This approach confines modelling to the smaller scales and can be naturally accommodated in a spectral element setting if one takes into account the hierarchical nature of the polynomial expansions used to represent the numerical solution \cite{munts2007modal,wasberg2009variational}. More recently, DG-based VMS approaches have received renewed attention, see \cite{collis2002dgvms,abba2015dynamic,chapelier2016development}. The effects of polynomial dealiasing on this kind of approach have been considered in \cite{beck2016influence}, which showed that under-integration can partially mask the physics embedded in the VMS model employed. An alternative strategy where dissipation is also confined to the small scales consists in the addition of spectral vanishing viscosity (SVV), see e.g.\ \cite{karamanos2000spectral,kirby2002coarse,kirby2006stabilisation}. SVV can be understood as a higher-order viscosity, designed to affect only the highest solution wavenumbers or polynomial modes \cite{pasquetti2005spectral,moura2016eigensolution}. This allows for the exponential convergence property of spectral element methods to be maintained, which helps achieving high accuracy at well-resolved laminar and transitional flow regions. Successful applications of this approach, more commonly used with CG methods, can be found in \cite{minguez2008high,koal2012adapting,lombard2015implicit,serson2016direct}. This SVV-based eddy-resolving approach can be considered an implicit LES (iLES) methodology, where, broadly speaking, numerical stabilization techniques are relied upon to dissipate small scales in lieu of a subgrid-scale turbulence model \cite{grinstein2007implicit}.

Implicit LES approaches based on discontinuous spectral element methods have also received significant attention in recent years \cite{parsani2010animplicit,uranga2011implicit,beck2014high,vermeire2016implicit}. This is partially because the upwind dissipation peculiar to these schemes already provides good stability for advection-dominated flows. As a result, valid solutions can be obtained without additional stabilization or modelling at moderate Reynolds numbers. The accuracy of these methods with regards to iLES is in great part due to a favourable (SVV-like) dissipative behaviour in wavenumber space, which does not affect the large scales directly and is only significant at high wavenumbers/frequencies \cite{van2007dispersion,van2008stability,vincent2011insights,moura2015,moura2016ontheeddy}. Still, as mentioned previously, this stabilization can be insufficient for high Reynolds number flows, especially at higher discretisation orders, even with over-integration. Polynomial filtering techniques can also be used for dealiasing. With the proper choice of the filter parameters the overall robustness is improved by suppressing spurious oscillations \cite{tcfd2012,flad2016simulation}.

In this work, we focus on issue (b) above. We consider consistent/over-integration as a reference dealiasing strategy, which in principle should only improve upon the accuracy and stability of DG-based iLES, and compare it against a novel class of split-form DG discretisations developed over recent years \cite{gassner_skew_burgers,gassner_kepdg,gassner2015,gassner2016split}. The latter are expected to have positive built-in dealiasing characteristics and have already demonstrated remarkable potential in terms of robustness. More specifically, in \cite{gassner2016split}, certain split forms proved capable of stabilizing the inviscid Taylor-Green vortex (TGV) problem \cite{shu2005numerical} even at very high polynomial orders. This test case is extremely demanding in terms of stability due to the absence of molecular viscosity and can easily diverge numerically regardless of consistent/over-integration \cite{moura2016ontheeddy}. However, the key goal beyond robustness is, of course, accuracy. The present study offers a detailed analysis of a comprehensive set of inviscid TGV test cases, and demonstrates that split form DG approaches can offer results of quality very similar to that achieved by consistent/over-integration. Although split forms are formulated by averaging standard conservative and advective forms of the governing equations, previous results derived in the context of finite difference methods \cite{fisher2013} can be invoked to guarantee that certain nodal split form DG discretisations remain conservative. From a numerical perspective, therefore, this approach is believed to aggregate significant advantages for eddy-resolving computations of high Reynolds number compressible flows. However, the main goal of this study is not to advocate DG-based iLES approaches. \comment{Instead the goal is, for the first time,} to investigate key aspects of solution quality and robustness the two dealiasing approaches \comment{for nodal DG methods} might exhibit when viscous effects are negligibly small.

This paper is organised as follows. Section~\ref{sec:DGSEM} presents a brief overview of a standard nodal DG formulation along with the two dealiasing strategies under consideration. Next, Section~\ref{sec:aliasingProperties} provides a preliminary analysis on the built-in dealiasing mechanisms of relevant split form DG methods in the context of Burgers' turbulence. In Section~\ref{sec:analysisTools} we introduce the inviscid TGV test cases and discuss the stability of the two dealiasing techniques considered. Section~\ref{sec:solQuality} directly compares the solution quality of the two approaches for various polynomial orders and different Riemann solvers. 
Concluding remarks are given in Section~\ref{sec:conclusions}.

\section{Discontinuous Galerkin discretisations}
\label{sec:DGSEM}

This section introduces the different DG formulations considered in this study. We first briefly review the standard and over-integrated DGSEM approaches in Section~\ref{sec:nodalDG}. We then discuss the split form DG discretisations which are relevant to our work in Section~\ref{sec:SplitForm}. Lastly, in Section~\ref{sec:numflux} we describe in more detail the numerical interface fluxes used with the different formulations being considered.

\subsection{Standard and over-integrated DGSEM approaches}
\label{sec:nodalDG}

We consider a system of equations given by conservation laws in the form
\begin{equation} \label{eq:ConsGeneral}
\frac{\partial \vec{q}}{\partial t} + \frac{\partial \vec{f}}{\partial x} + \frac{\partial \vec{g}}{\partial y} + \frac{\partial \vec{h}}{\partial z} = \vec{0} ,
\end{equation}
where $\vec{q} (x,y,z,t)$ is the vector of conserved variables and $\vec{f}(\vec{q})$, $\vec{g}(\vec{q})$ and $\vec{h}(\vec{q})$ are the flux vectors governing the transport of $\vec{q}$ in a physical domain $\Omega$. For the spectral element approximation, we subdivide $\Omega$ into non-overlapping elements $\Omega_e$ such that $\bigcup_e \Omega_e$ is a mesh of the computational domain. For simplicity, we restrict the following presentation to straight-sided hexahedral elements, but remark that the extension to more general curvilinear meshes is available, see e.g.\ \cite{Kopriva:2009nx}.

In the nodal DGSEM approach considered, the solution is approximated within each element $\Omega_e$ through a polynomial expansion such as
\begin{equation} \label{eq:polyApprox}
\vec{q} \approx \sum_{n} \check{\vec{q}}_{n} (t) \, \varphi_{n} (x,y,z) ,
\end{equation}
in which the summation index is assumed to span a polynomial space of degree $N$, and where time-dependent coefficients $\check{\vec{q}}_{n} (t)$ are associated to the nodal basis functions $\varphi_{n} (\vec{x})$. The spatial dependence usually relies on a simple trilinear mapping relation $\vec{x} = \vec{x} ( \vec{\xi} )$, cf.\ \cite{gassner2016split}, that connects the physical space $\Omega_e$ to the standard domain $\Omega_0 = [-1, 1]^3$. This allows for the basis functions to be defined with regards to $\Omega_0$ directly, namely
\begin{equation}
\varphi_{n} ( \vec{x} (\vec{\xi}) ) = \ell_i (\xi) \ell_j (\eta) \ell_k (\zeta) ,\quad n = 1,\ldots,(N+1)^3,
\end{equation}
where, for example, index $i = i(n)$ references the $N$-th order Lagrange polynomial $\ell_i (\xi)$,
\begin{equation}\label{eq:lagrange_basis}
\ell_i(\xi) = \prod\limits_{\underset{l \neq i}{l = 0}}^N \frac{\xi - \xi_l}{\xi_i - \xi_l} , \qquad \mbox{with } i = 0, \dots, N,
\end{equation}
which is based on an appropriate set of nodes $\{ \xi_l \}^N_{l=0}$ defined within $[-1, 1]$. The nodal DG methods in this work uses a collocation approach, e.g. \cite{Kopriva:2009nx}, where interpolation and quadrature nodes are selected to be the Legendre-Gauss-Lobatto (LGL) nodes. The choice of LGL quadrature is important, because it ensures that the discrete DG derivative matrix and the discrete mass
matrix satisfies the summation-by-parts (SBP) property for any polynomial order \cite{gassner_skew_burgers}. The SBP property is crucial to define a split form DGSEM that remains conservative \cite{fisher2013_2,fisher2013}. Finally, we note the cardinal property of Lagrange polynomials \eqref{eq:lagrange_basis}, whereby $\ell_i(\xi_l) = \delta_{il}$, with $\delta_{il} = 1$ for $i = l$ and $\delta_{il} = 0$ otherwise.

We then require that (\ref{eq:ConsGeneral}), when evaluated with (\ref{eq:polyApprox}), should vanish locally in a projection sense, i.e.\
\begin{equation} \label{eq:weakFormElement}
\int_{\Omega_e} ( \vec{q}_t + \vec{f}_x + \vec{g}_y + \vec{h}_z ) \, \varphi_m \, \mathrm{d}\vec{x} = \vec{0} ,\quad m =1,\ldots,(N+1)^3,
\end{equation}
for all the basis functions, $\varphi_m$, used in (\ref{eq:polyApprox}). Having in mind the simulations to be discussed in this study, we now restrict our discretisation to the case of Cartesian meshes, where element-wise Jacobian factors and metric terms simplify considerably to $J = \frac{1}{8} \Delta x \Delta y \Delta z$, $x_{\xi} = \half \Delta x$, $y_{\eta} = \half \Delta y$, $z_{\zeta} = \half \Delta z$, in which $\Delta x$, $\Delta y$, and $\Delta z$ represent the element side lengths along the three Cartesian directions. In this case, casting (\ref{eq:weakFormElement}) as an integration over $\Omega_0$ yields \cite{Kopriva:2009nx}
\begin{equation} \label{eq:weakFormElementRef}
\int_{\Omega_0} ( J \vec{q}_t + \nabla_{\xi} \cdot \mathscr{F} ) \, \varphi_m \, \mathrm{d}\vec{\xi} = \vec{0} ,\quad m =1,\ldots,(N+1)^3,
\end{equation}
where $\mathscr{F} = [ \mathcal{F}, \mathcal{G}, \mathcal{H} ]^T$ is a vector containing the contravariant fluxes, which incorporate the element-wise constant metric terms, namely
\begin{equation} \label{eq:newVars}
\mathcal{F} (\vec{q}) = y_\eta z_\zeta \, \vec{f} (\vec{q}) , \quad \mathcal{G} (\vec{q}) = x_\xi z_\zeta \, \vec{g} (\vec{q}) , \quad \mathcal{H} (\vec{q}) = x_\xi y_\eta \, \vec{h} (\vec{q}) .
\end{equation}
Finally, we integrate (\ref{eq:weakFormElementRef}) by parts and replace the boundary terms stemming from the divergence theorem with numerical interface fluxes, $\mathscr{F}^*$. The resulting statement for the weak DG formulation reads
\begin{equation} \label{eq:IBPWeakForm}
J \int_{\Omega_0} \vec{q}_t \, \varphi_m \, \mathrm{d}\boldsymbol{\xi} + \int_{\partial \Omega_0} \varphi_m \, \mathscr{F}^* \cdot \mathrm{d}\vec{S} - \int_{\Omega_0} \mathscr{F} \cdot \nabla_{\xi} \, \varphi_m \, \mathrm{d}\boldsymbol{\xi} = \vec{0} ,\quad m =1,\ldots,(N+1)^3,
\end{equation}
in which the differential surface vector $d\vec{S}$ points toward the outside of $\Omega_0$. The DG schemes rely on well-known Riemann solvers \cite{torobook} for $\mathscr{F}^*$ to resolve discontinuities in the approximation at element surfaces. There is additional flexibility for a split form DGSEM which can use special numerical fluxes to recover certain secondary properties like kinetic energy preservation \cite{gassner2016split}, as will be discussed in Section~\ref{sec:numflux}. 

Each integral in (\ref{eq:IBPWeakForm}) is evaluated numerically with a Gauss type quadrature, given the number $Q$ of integration nodes, see e.g.\ \cite{KarniadakisSherwin}. The \textit{standard} DGSEM relies on collocation with $Q = N+1$ integration points per direction, usually LGL nodes \cite{Kopriva:2009nx,hesthaven2007}. While this choice only guarantees exact integration of flux functions depending linearly on the solution, the same choice is adopted, in many cases, for problems involving non-linear flux functions. This collocated approximation is advantageous in terms of computational cost, cf.\ \cite{Hindenlang2012}, but can also be numerically unstable and less accurate due to polynomial aliasing errors introduced by the under-integration of non-linear flux functions, especially for under-resolved computations \cite{Kirby2003,tcfd2012,mengaldo2015}. The straightforward way to alleviate this problem is to simply use a larger number of quadrature points, i.e.\ $Q > N+1$, in what would be an \textit{over-integrated} DGSEM.

The over-integration mentioned above becomes a consistent integration when the quadrature order is chosen to be consistent with the nonlinearity of flux functions. For example, one typically requires $Q \approx \frac{3}{2} (N+1)$ for quadratic and $Q \approx 2 (N+1)$ for cubic nonlinearities. These are respectively common choices for the consistent integration (of the advective terms) of the incompressible and compressible Navier-Stokes equations \cite{KarniadakisSherwin}. For compressible flows with small density variations, using $Q \approx \frac{3}{2} (N+1)$ might also reproduce the effects of consistent integration \cite{beck2014high,beck2016influence}. However, we remark that Gauss type quadratures can only guarantee exact results for the integration of polynomial functions, and that compressible flow equations (written in the usual conservative variables) have flux vectors whose components are actually rational functions. As a result, strictly speaking, the associated quadratures are not exact even for $Q \approx 2 (N+1)$. As asymptotic convergence to exact results is nevertheless guaranteed \cite{trefethen2013approximation}, quadrature errors can be expected to decay numerically to negligible levels for $Q$ sufficiently large, as exemplified in Section~\ref{sec:analysisTools}.

The standard DGSEM is perhaps more commonly known in the so-called strong form variant of (\ref{eq:IBPWeakForm}), which is obtained with an additional integration by parts and reads
\begin{equation} \label{eq:IBPAgain}
J \int_{\Omega_0} \vec{q}_t \, \varphi_m \, \mathrm{d}\boldsymbol{\xi} + \int_{\partial \Omega_0} \varphi_m \, \left( \mathscr{F}^* - \mathscr{F} \right) \cdot \mathrm{d}\vec{S} + \int_{\Omega_0} \varphi_m \, \nabla_{\xi} \cdot \mathscr{F} \, \mathrm{d}\boldsymbol{\xi} = \vec{0} ,\quad m=1,\ldots,(N+1)^3,
\end{equation}
which in general differs from (\ref{eq:IBPWeakForm}) when quadratures are inexact. For standard DGSEM, however, it is possible to show that forms (\ref{eq:IBPWeakForm}) and (\ref{eq:IBPAgain}) are actually discretely equivalent \cite{KoprivaGassner_GaussLob}, despite any quadrature errors. Also, due to the cardinal property of the Lagrange polynomials employed and the collocation of quadrature and solution nodes, the mass matrix stemming from the leftmost integral above has a diagonal structure and is given by \cite{hesthaven2007,Kopriva:2009nx}
\begin{equation} \label{eq:massMat}
\mathbf{M} = \textrm{diag} (\omega_0, \dots, \omega_N) ,
\end{equation}
where $\{ \omega_l \}^N_{l=0}$ is the set of quadrature weights associated to the set of LGL nodes $\{ \xi_l \}^N_{l=0}$ adopted. We note that when consistent integration is employed, mass matrices of nodal approaches are actually full. This is what causes dispersion--diffusion characteristics of standard and over-integrated DGSEM approaches to differ at higher wavenumbers \cite{hesthaven2007}.

In the case of standard DGSEM, the diagonal mass matrix allow for the semi-discrete evolution equation of each nodal coefficient to be written independently, in a form that resembles (\ref{eq:ConsGeneral}), namely
\begin{equation} \label{eq:nodalDGSEM}
J( \check{\vec{q}}_t )_{ijk} + ( \widetilde{\mathcal{F}}_\xi )_{ijk} + ( \widetilde{\mathcal{G}}_\eta )_{ijk} + ( \widetilde{\mathcal{H}}_\zeta )_{ijk} = \vec{0} ,
\end{equation}
in which the modified (tilde) fluxes above incorporate the interface contributions and are given by
\begin{align}
(\widetilde{\mathcal{F}}_\xi)_{ijk} &= \frac{1}{\omega_i} \left[ \delta_{iN} \left( \mathcal{F}^* - \mathcal{F} \right)_{Njk} - \delta_{0i} \left( \mathcal{F}^* - \mathcal{F} \right)_{0jk} \right] + \sum_{n=0}^N \mathbf{D}_{in} (\mathcal{F})_{njk} , \label{eq:RHSStrong1} \\
(\widetilde{\mathcal{G}}_\eta)_{ijk} &= \frac{1}{\omega_j} \left[ \delta_{jN} \left( \mathcal{G}^* - \, \mathcal{G} \right)_{Njk} - \delta_{0j} \left( \mathcal{G}^* - \, \mathcal{G} \right)_{0jk} \right] \, + \, \sum_{n=0}^N \mathbf{D}_{jn} (\mathcal{G})_{ink} , \label{eq:RHSStrong2} \\
(\widetilde{\mathcal{H}}_\zeta)_{ijk} &= \frac{1}{\omega_k} \left[ \delta_{kN} \left( \mathcal{H}^* - \mathcal{H} \right)_{ijN} - \delta_{0k} \left( \mathcal{H}^* - \mathcal{H} \right)_{ij0} \right] + \sum_{n=0}^N \mathbf{D}_{kn} (\mathcal{H})_{ijn} , \label{eq:RHSStrong3}
\end{align}
where $\mathbf{D}$ is the standard polynomial derivative matrix, see e.g.\ \cite{Kopriva:2009nx}, defined as
\begin{equation} \label{eq:standardDer}
\mathbf{D}_{mn} = \frac{\partial \ell_n}{\partial \xi} \Bigg|_{\xi = \xi_m} , \qquad \mbox{for } m, n = 0, \dots, N.
\end{equation}
In the following section, different split form approaches will be derived based on simple modifications of the rightmost terms in (\ref{eq:RHSStrong1})-(\ref{eq:RHSStrong3}).

\subsection{Split form DG discretisations}
\label{sec:SplitForm}

Here we discuss the so-called split form approaches to the discretisation of the advective terms governing general fluid flow. The structure of any split formulations is equation dependent because knowledge of the type of non-linearity (e.g. integer powers of the fluid velocity) is necessary to create an average of the conservative and advective forms of the advective terms \cite{ducros2000,sjogreen2010,gassner2016split}. In this work we consider the compressible Euler equations, most commonly stated in conservation form as (\ref{eq:ConsGeneral}) with $\vec{q} = [ \rho, \rho u, \rho v, \rho w, \rho e ]^T$ and
\begin{equation} \label{eq:compr-fluxes}
\vec{f} (\vec{q}) = \! \begin{bmatrix}
\rho\,u\\
\rho\,u^2 +p\\
\rho\,u\,v\\
\rho\,u\,w\\
(\rho\,e+p) \, u
\end{bmatrix} , \quad
\vec{g} (\vec{q}) = \! \begin{bmatrix}
\rho\,v\\
\rho\,u\,v \\
\rho\,v^2+p\\
\rho\,v\,w\\
(\rho\,e+p) \, v
\end{bmatrix} , \quad
\vec{h} (\vec{q}) = \! \begin{bmatrix}
\rho\,w\\
\rho\,u\,w\\
\rho\,v\,w\\
\rho\,w^2+p\\
(\rho\,e+p) \, w
\end{bmatrix} ,
\end{equation}
where $\rho\,e = p/(\gamma-1) + \rho\,(u^2+v^2+w^2)/2$ while $\rho$, $u$, $v$, $w$, $p$, $e$ and $\gamma = 7/5$ stand respectively for density, the three velocity components along the Cartesian directions, pressure, specific total energy and the usual ratio of specific heats.

Split formulations are alternative discretisations of the non-linear transport terms typical of many partial differential equations (PDEs). Essentially, splitting techniques are built by averaging conservative and advective forms of a PDE. Relevant split forms usually improve upon the robustness of numerical methods by reducing aliasing errors. Split formulations have been often used in conjunction with high-order finite difference schemes \cite{FDaliasing,ducros2000,larsson2007,kennedy2008,Morinishi2010276} and spectral methods \cite{blaisdell1996effect,Zang199127}. There are many ways to rewrite the non-linear advective terms of a PDE. A good overview of different split form approaches can be found in \cite{pirozzoli2011}. Depending on how one interprets the non-linearity of the Euler fluxes (quadratic, cubic or rational) there are several ways to rewrite the equations in an equivalent split form. To simplify the discussion, we introduce notation for specific instances of the quadratic splitting 
\begin{equation} \label{eq:quadSplit}
(a \cdot b)_x \coloneqq \frac{1}{2} \, (a \, b)_x + \frac{1}{2} \left( a_x \, b + a \, b_x \right) ,
\end{equation}
and of the cubic splitting explored by Kennedy and Gruber \cite{kennedy2008}
\begin{equation} \label{eq:cubicSplit}
(a \cdot b \cdot c)_x \coloneqq \frac{1}{4} \, (a \, b \, c)_x + \frac{1}{4} \left[ a_x \, (b \, c) + b_x \, (a \, c) + c_x \, (a \, b) \right] + \frac{1}{4} \left[ a \, (b \, c)_x + b \, (a \, c)_x + c \, (a \, b)_x \right] .
\end{equation}

Some splitting strategies are motivated by the type of flow physics that the numerics are expected to capture. For example, Ducros et al.\ \cite{ducros2000} considered flows with small density variations and used the quadratic splitting \eqref{eq:quadSplit} for each term in the Euler fluxes. Alternatively, Kennedy and Gruber \cite{kennedy2008} wanted to account for density variations and hence used both \eqref{eq:quadSplit} and \eqref{eq:cubicSplit}. These two splitting strategies, hereafter denoted respectively as DU and KG, are given, e.g.\ in the $x$-direction, by
\begin{equation} \label{eq:op_DU}
\vec{f}_x^{DU} (\vec{q}) =
\begin{bmatrix}
\half \left[ (\rho u)_x + \rho(u)_x + u(\rho)_x \right] \\[0.1cm]
\frac{1}{4} \left[ (\rho u^2)_x + \rho u (u)_x + u(\rho u)_x \right] + p_x \\[0.1cm]
\frac{1}{4} \left[ (\rho v u)_x + \rho v (u)_x + u(\rho v)_x \right] \\[0.1cm]
\frac{1}{4} \left[ (\rho w u)_x + \rho w (u)_x + u(\rho w)_x \right] \\[0.1cm]
\frac{1}{4} \left[ (\rho e u + p u)_x + (\rho e + p) (u)_x + u(\rho e + p)_x \right]
\end{bmatrix},
\end{equation}
and
\begin{equation} \label{eq:op_KG}
\vec{f}_x^{KG} (\vec{q}) = \begin{bmatrix}
\half \left[ (\rho u)_x + \rho(u)_x + u(\rho)_x \right] \\[0.1cm]
\frac{1}{4} \left[ (\rho u^2)_x + \rho (u^2)_x + 2u(\rho u)_x + u^2(\rho)_x + 2\rho u(u)_x \right] + p_x \\[0.1cm]
\frac{1}{4} \left[ (\rho u v)_x + \rho (uv)_x + u(\rho v)_x + v(\rho u)_x + uv(\rho)_x + \rho v(u)_x + \rho u(v_x) \right] \\[0.1cm]
\frac{1}{4} \left[ (\rho u w)_x + \rho (uw)_x + u(\rho w)_x + w(\rho u)_x + uw(\rho)_x + \rho w(u)_x + \rho u(w_x) \right] \\[0.1cm]
\frac{1}{4} \left[ (\rho e u)_x + \rho(e u)_x + e(\rho u)_x + u(\rho e)_x+ eu(\rho)_x + \rho u (e)_x + \rho e (u)_x \right] + \half \left[ (pu)_x + p(u)_x + u(p_x) \right]	\end{bmatrix},
\end{equation}
or, by using the compact notation introduced in (\ref{eq:quadSplit}) and (\ref{eq:cubicSplit}),
\begin{equation} \label{eq:compact_DUKG}
\vec{f}_x^{DU} (\vec{q}) = \! \begin{bmatrix}
\rho \cdot u\\
\rho u \cdot u + p\\
\rho v \cdot u\\
\rho w \cdot u\\
(\rho e + p) \cdot u
\end{bmatrix}_x , \quad
\vec{f}_x^{KG} (\vec{q}) = \! \begin{bmatrix}
\rho \cdot \, u\\
\rho \cdot u \cdot u + p\\
\rho \cdot u \cdot v\\
\rho \cdot u \cdot w\\
\rho \cdot e \cdot u + p \cdot u
\end{bmatrix}_x .
\end{equation}
We note that while the DU splitting does not lead to discretisations that are formally kinetic energy preserving \cite{pirozzoli2011}, the KG splitting does allow for that possibility \cite{jameson2008}.

With the concept of split forms, it is possible to build alternative DGSEM discretisations with special properties guaranteed at the discrete level, such as kinetic energy preservation or entropy consistency \cite{gassner_skew_burgers,gassner_kepdg,gassner2015,gassner2016split}. However, because these are constructed by averaging conservative and non-conservative forms of the governing equations, it is important that the resulting DG scheme remains locally conservative. This can be guaranteed for the forms considered here by the following. We start by modifying the computation of the volume integral terms of the standard DGSEM in its strong variant, cf.\ \eqref{eq:IBPAgain} and (\ref{eq:RHSStrong1})-(\ref{eq:RHSStrong3}). The collocated mass \eqref{eq:massMat} and derivative \eqref{eq:standardDer} matrices are combined into the matrix $\qmat \coloneqq \mmat \dmat$, which has the summation-by-parts (SBP) property \cite{carpenter_esdg}, namely
\begin{equation}
\qmat + \qmat^T = \mat{B} \coloneqq \textrm{diag} (-1, 0, \dots, 0, 1),
\end{equation}
that is used to mimic integration-by-parts at a discrete level, by manipulating the derivative matrix as
\begin{equation}
\dmat = \mmat^{-1} \qmat = \mmat^{-1} \mat{B} - \mmat^{-1} \qmat^T .
\end{equation}
A remarkable result presented in Fisher et al. \cite{fisher2013_2} and Fisher and Carpenter \cite{fisher2013} showed that diagonal norm SBP matrices such as $\dmat$ can be reinterpreted as sub-cell finite volume type differencing operators. This result is what guarantees local conservation of diagonal norm SBP discretisations in the sense of Lax-Wendroff due to a telescoping flux differencing \cite{fisher2013}. For more details and proofs, the interested reader is referred to \cite{gassner2016split}.

As further explored in \cite{gassner2016split}, the flux differencing formulation mentioned above can be recast into a matrix type product form by the introduction of auxiliary sub-cell numerical flux functions $\mathcal{F}^{\#}$, $\mathcal{G}^{\#}$, and $\mathcal{H}^{\#}$, which are required to be symmetric and consistent, e.g.\
\begin{equation}
\mathcal{F}^{\#} ( \vec{q}_L, \vec{q}_R ) = \mathcal{F}^{\#} ( \vec{q}_R, \vec{q}_L ), \quad \mbox{ and } \quad \mathcal{F}^{\#} ( \vec{q}, \vec{q} ) = \mathcal{F} ( \vec{q} ) .
\end{equation}
Finally, this matrix type product form based on sub-cell numerical fluxes can be linked (via transitivity) to the original element-wise differentiation operation, cf.\ (\ref{eq:RHSStrong1})-(\ref{eq:RHSStrong3}), such that
\begin{align}
&\sum_{n=0}^N \mathbf{D}_{in} (\mathcal{F})_{njk} \, \approx \, \,2\, \sum_{n=0}^N \mathbf{D}_{in} \, \mathcal{F}^{\#} ( \vec{q}_{ijk}, \vec{q}_{njk} ) , \label{eq:splitStrong1} \\
&\sum_{n=0}^N \mathbf{D}_{jn} (\mathcal{G})_{ink} \,\, \approx \, \,2\, \sum_{n=0}^N \mathbf{D}_{jn} \, \mathcal{G}^{\#} ( \vec{q}_{ijk}, \vec{q}_{ink} ) , \label{eq:splitStrong2} \\
&\sum_{n=0}^N \mathbf{D}_{kn} (\mathcal{H})_{ijn} \, \approx \, \,2\, \sum_{n=0}^N \mathbf{D}_{kn} \, \mathcal{H}^{\#} ( \vec{q}_{ijk}, \vec{q}_{ijn} ) , \label{eq:splitStrong3}
\end{align}
where the original differentiation is recovered for fully central sub-cell numerical fluxes, e.g.\ $\mathcal{F}^{\#} ( \vec{q}_{ijk}, \vec{q}_{njk} ) = \half [ \mathcal{F} ( \vec{q}_{ijk} ) + \mathcal{F} ( \vec{q}_{njk} )]$, but also different split forms can be implemented effortlessly with other choices. For instance, if we use instead the product of two averages or the product of three averages for the sub-cell numerical flux we recover discrete versions of the quadratic \eqref{eq:quadSplit} and cubic \eqref{eq:cubicSplit} split forms, see \cite{gassner2016split} for complete details. This creates a dictionary where one can easily translate between the continuous split forms and their discrete counterparts. For example, in the $x$-direction, the split forms in \eqref{eq:compact_DUKG} can be recovered respectively through the sub-cell numerical fluxes $\mathcal{F}^{\#}_{DU} ( \vec{q}_{ijk}, \vec{q}_{njk} )$ and $\mathcal{F}^{\#}_{KG} ( \vec{q}_{ijk}, \vec{q}_{njk} )$ given by
\begin{equation} \label{eq:discrete_DUKG}
\mathcal{F}^{\#}_{DU} = y_{\eta} z_{\zeta}
\begin{bmatrix}
\avg{\rho} \avg{u} \\[0.1cm]
\avg{\rho u} \avg{u} + \avg{p} \\[0.1cm]
\avg{\rho v} \avg{u} \\[0.1cm]
\avg{\rho w} \avg{u} \\[0.1cm]
\avg{\rho e + p} \avg{u}
\end{bmatrix} ,
\qquad
\mathcal{F}^{\#}_{KG} = y_{\eta} z_{\zeta}
\begin{bmatrix}
\avg{\rho} \avg{u} \\[0.1cm]
\avg{\rho} \avg{u}^2 + \avg{p} \\[0.1cm]
\avg{\rho} \avg{u} \avg{v} \\[0.1cm]
\avg{\rho} \avg{u} \avg{w} \\[0.1cm]
\avg{\rho} \avg{e} \avg{u} + \avg{p} \avg{u}
\end{bmatrix} ,
\end{equation}
where $\avg{\cdot}$ represents the arithmetic mean --- in this case, between the relevant quantities from $\vec{q}_{ijk}$ and $\vec{q}_{njk}$. The dot splitting notation (\ref{eq:quadSplit}) and(\ref{eq:cubicSplit}) elucidates the usage of symmetric sub-cell fluxes: one simply replaces the dots with products of arithmetic means. Reference \cite{gassner2016split} discusses many other sub-cell flux functions, e.g.\ those related to other well-known split forms like that of Morinishi \cite{Morinishi2010276} or Pirozzoli \cite{pirozzoli2011}.

Lastly, we replace the rightmost terms in (\ref{eq:RHSStrong1})-(\ref{eq:RHSStrong3}) with (\ref{eq:splitStrong1})-(\ref{eq:splitStrong3}) and, following \cite{gassner2016split}, we connect the choice of the sub-cell flux to that of the interface numerical flux. For example, if the KG splitting is adopted, the local Lax-Friedrichs (LLF) interface flux function will have, e.g.\ for the $x$-direction, the adapted form
\begin{equation} \label{eq:dissipationTerm}
\mathcal{F}_{LLF}^* (\vec{q}_L , \vec{q}_R) \, = \, \mathcal{F}_{KG}^{\#} (\vec{q}_L , \vec{q}_R) \, - \, \half \, \lambda_{max} \, ( \vec{q}_R - \vec{q}_L ) ,
\end{equation}
where $\lambda_{max}$ is an estimate of the fastest wave speed (in absolute value) at the interface between the left and right solution states $\vec{q}_L$ and $\vec{q}_R$. More details about the interface flux functions employed in the different DGSEM approaches under consideration will be discussed in the next section.

\subsection{Selection of the numerical interface flux function}
\label{sec:numflux}

We briefly outline the two types of interface flux functions that will be used with the over-integrated and split form DGSEM approaches considered. For over-integrated DGSEM, we use the classical forms of LLF and Roe-type schemes \cite{torobook}. For example, in the $x$-direction, the LLF flux formula reads
\begin{equation} \label{eq:LLFOI}
\mathcal{F}_{LLF}^* (\vec{q}_L , \vec{q}_R) \, = \, \half \, [ \mathcal{F} (\vec{q}_L) + \mathcal{F} (\vec{q}_R) ] \, - \, \half \, \lambda_{max} \, ( \vec{q}_R - \vec{q}_L ) ,
\end{equation}
where terms are denoted analogously to \eqref{eq:dissipationTerm}. As for the classical Roe solver, the dissipative (rightmost) term relies on a characteristic decomposition of the flux Jacobian evaluated at Roe's average state \cite{Roe:1981:JCP81}.

In contrast, for split form DGSEM, we couple the choices of sub-cell and interface numerical fluxes, as shown for the LLF case in \eqref{eq:dissipationTerm}. When the Roe flux is employed with the DU splitting, the symmetric part of Roe's flux formula is replaced with $\mathcal{F}_{DU}^{\#} (\vec{q}_L , \vec{q}_R)$ and the original dissipation term is maintained. In case the KG splitting is adopted, however, a slight modification is made to the dissipation term. We note that the specific Roe averaging remains identical to the classical scheme, but the dissipation term is modified to ensure that the KG scheme is kinetic energy stable, see \cite{jameson2008,chandrashekar2013,gassner2016split}. To do so, we alter the diagonal matrix of eigenvalues so that the first eigenvalue matches the last one \cite{chandrashekar2013}, namely
\begin{equation} \label{eq:newDiag}
\mathbf{\Lambda} = \textrm{diag} ([ \, u+a, u, u, u, u+a \, ]) ,
\end{equation}
where $a$ represents the speed of sound.

Finally, for the split form schemes we note a somewhat remarkable property. It is possible, in some cases, to obtain numerical solutions of the compressible Euler equations which nearly conserve the total kinetic energy \cite{gassner2016split}. This can be achieved, for example, with a ``fully central'' interface flux such as
\begin{equation} \label{eq:centralKG}
\mathcal{F}^*_{central} \, = \, \mathcal{F}^\#_{KG} ,
\end{equation}
which lacks the usual interface dissipation term. Simulations based on this central flux will be discussed in Section\ \ref{sec:strangeQuality}, where solutions obtained with very high polynomial orders will be compared to lower order solutions that conserve total kinetic energy.

\comment{All simulation results in Secs. \ref{sec:analysisTools} and \ref{sec:solQuality} are integrated in time with an explicit five stage, fourth order accurate low storage Runge-Kutta scheme \cite{Carpenter&Kennedy:1994}, where a stable time step is computed according to an adjustable coefficient CFL$\in(0,1]$, the local maximum wave speed and the relative grid size, e.g. \cite{gassner2011}.}

\section{Analysis of polynomial aliasing errors of split form approaches}
\label{sec:aliasingProperties}

We now undertake a preliminary assessment of how inexact quadratures can affect the energy of different solution coefficients and how the split form approach tends to keep this energy ``under control," thus increasing numerical robustness. The objective is to explain, at least partially, why certain split form DG approximations have positive stabilisation and dealiasing properties. A similar analysis has been carried out for pseudo-spectral methods (in Fourier space) by Blaisdell et al. in \cite{blaisdell1996effect}, where it was demonstrated that split form discretisations suppress aliasing errors as if they had a built-in dealiasing mechanism. Here we simply illustrate through numerical experiments that this convenient feature seems to hold for nodal DG-based spectral element methods. More specifically, we consider the Burgers' equation in one spatial dimension, given by \eqref{eq:ConsGeneral} with flux function $f(q) = q^2/2$. Similar numerical experiments have been used in previous works \cite{Kirby2003,kirby2006stabilisation} to illustrate the effect of different dealiasing techniques in controlling the solution's modal energy.

The inner-product formulation of the strong split form version of the Burgers' equation reads \cite{gassner_skew_burgers}
\begin{equation} \label{eq:split_quads_strong}
	\left\langle \frac{\partial q}{\partial t} , \ell_i \right\rangle_Q = -\, \alpha \left\langle \frac{(q^2)'}{2} , \ell_i \right\rangle_Q +\, (\alpha - 1) \left\langle q q' , \ell_i \right\rangle_Q - \left[ \left(f^* - \frac{q^2}{2}\right) \ell_i \right]^{\oplus}_{\ominus} , \,\,\,\,\,\, \mbox{for $i = 0, \dots, N$,}
\end{equation}
where $(\cdot)'$ denotes differentiation in space and $\langle \cdot , \cdot \rangle_Q$ stands for the quadrature of the product between two functions using $Q$ integration points (LGL nodes). Furthermore, test functions $\ell_i$ are the usual Lagrange polynomials of degree $N$, cf.\ \eqref{eq:lagrange_basis}, and $f^*$ is the numerical flux required at the right and left elemental interfaces, denoted above as $\oplus$ and $\ominus$, respectively. The parameter $\alpha \in [0,1]$ is used to introduce biased averages between conservative ($\alpha = 1$) and non-conservative ($\alpha = 0$) discretisations.

In order to focus on the terms requiring quadratures, we work instead with the weak split form version of the Burgers' equation, where boundary contributions can be more easily dismissed. In the context of DGSEM discretisations, this version is entirely equivalent to (\ref{eq:split_quads_strong}), see \cite{gassner_skew_burgers}. If we disregard boundary contributions by setting $f^* = 0$, this weak form becomes
\begin{equation} \label{eq:split_quads}
\left\langle \frac{\partial q}{\partial t} , \ell_i \right\rangle_Q = \, \alpha \left\langle \frac{q^2}{2} , \ell'_i \right\rangle_Q +\, (\alpha - 1) \left[ \left\langle \frac{q q'}{2} , \ell_i \right\rangle_Q - \left\langle q , \frac{(q \ell_i)'}{2} \right\rangle_Q \right] - \cancel{ \left( f^* \ell_i \right)_{\ominus}^{\oplus} } .
\end{equation}

\comment{Note that in both formulations above, numerical differentiation of quadratic terms like $(q^2)'$ and $(q \ell_i)'$ is performed in accordance with the collocated nodal DGSEM framework: arguments are evaluated via point-wise product at the $Q$ integration points and treated as Lagrange polynomials of degree $N = Q-1$ (instead of $2N$), whose derivative is then taken in a standard way.}

\comment{Although our focus is on the characteristics of collocated nodal DGSEM,} the assessment of solution coefficients in modal space is particularly insightful as it allows for interpretation by analogy with Fourier space, although it is recognized that this analogy is somewhat limited. This will require the appropriate transformations between nodal and modal spaces to be employed. In the following, we restrict ourselves to the reference domain $[-1, 1]$ and consider a ``frozen'' numerical solution defined by the modal expansion
\begin{equation} \label{eq:modal_expansion}
	q^\delta (\xi) = \sum_{j=0}^{N} { \hat{q}_j L_j(\xi) } \,\,\, \mbox{ with } \,\,\, \hat{q}_j = (j+1)^{-5/6} \mbox{,}
\end{equation}
where $L_j(\xi)$ is the orthonormal Legendre polynomial of degree $j$. By the Fourier/modal space analogy, the above expansion is representative of a turbulent velocity field, since $\hat{q}_j^2 \propto (j+1)^{-5/3}$ as typical of turbulence energy spectra. We note that such scaling is achievable even with the Burgers' equation by means of a specific forcing \cite{zikanov1997,moura2015}. The coefficients of a nodal expansion equivalent to that of \eqref{eq:modal_expansion} are given by $\check{q} = \mathbf{V} \hat{q}$, where $\mathbf{V}_{ij} = L_j(\xi_i)$ is the Vandermonde matrix evaluated at the ($N+1$) LGL nodes $\xi_i$.

The quadratures in \eqref{eq:split_quads} are evaluated with the same set of $Q=N+1$ points $\xi_i$ mentioned above, in accord with the collocated nodal DGSEM approach, as follows:
\comment{\begin{equation}
	\langle q^2 , \ell'_i \rangle_{N+1} \, \approx \, \sum_{k=0}^{N} {q^\delta (\xi_k) \, q^\delta (\xi_k) \, \partial_{\xi}\ell_i (\xi_k) \, \omega_k} \mbox{,}
\end{equation}
\begin{equation}
	\langle q \, q' , \ell_i \rangle_{N+1} \, \approx \, \sum_{k=0}^{N} {q^\delta (\xi_k) \, \partial_{\xi}q^\delta (\xi_k) \, \ell_i (\xi_k) \, \omega_k} \mbox{,}
\end{equation}
\begin{equation}
	\langle (q \, \ell_i)' , q \rangle_{N+1} \, \approx \, \sum_{k=0}^{N} {\partial_{\xi}I^{\mathtt{LGL}} (\xi_k) \, q^\delta (\xi_k) \, \omega_k} \mbox{,}
\end{equation}}
where $\partial_{\xi} (\cdot) = \partial (\cdot) / \partial \xi$, $\omega_k$ are the quadrature weights and $I^{\mathtt{LGL}}$ is the LGL interpolant of the product $q^\delta \, \ell_i$, i.e.\ $I^{\mathtt{LGL}} (\xi)$ is a Lagrange polynomial of degree $N$ defined as $q^\delta (\xi_k) \, \ell_i (\xi_k)$ at the $N+1$ quadrature nodes. Note that $q^\delta (\xi) \, \ell_i (\xi)$ is actually a polynomial of degree $2N$, and so the approximation associated with the interpolant $I^{\mathtt{LGL}}$ is presumably related to the positive dealiasing properties observed in relevant split form DGSEM discretisations.

\begin{figure}[!ht]
	\centering
	\includegraphics[keepaspectratio=true, trim={0cm 0cm 0cm 0cm}, clip, width=16.5cm]{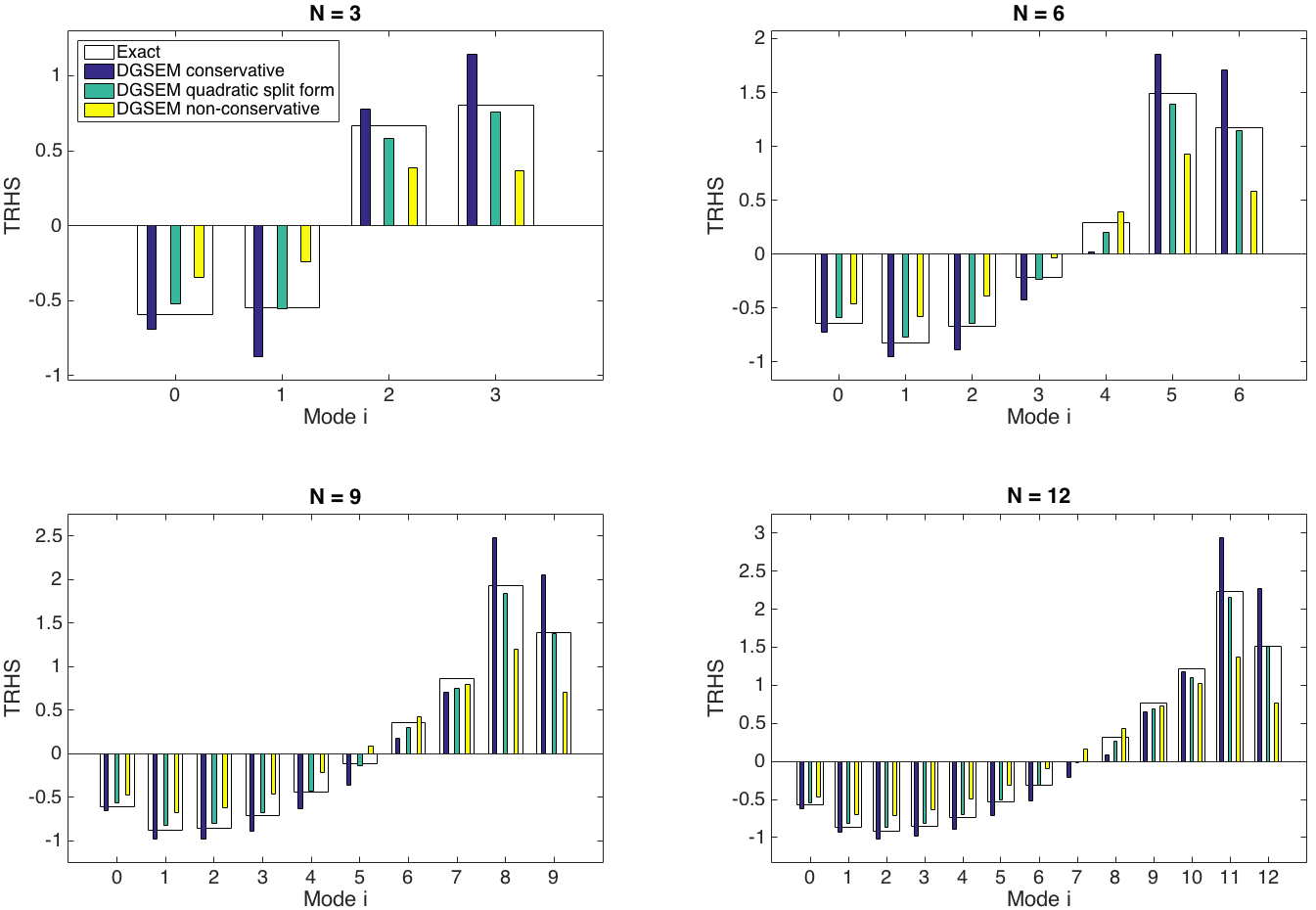}
	\caption{Comparison of integration errors for the 1D Burgers' equation with turbulent-like solution. The transformed (modal space) right-hand side of \eqref{eq:split_quads} is shown vs.\ the polynomial mode $i$ for different polynomial orders $N$. Conservative, quadratic split form and non-conservative DGSEM forms are compared against reference values obtained with exact integration.}
	\label{fig:1dComparison}
\end{figure}

The split form obtained from \eqref{eq:split_quads} with $\alpha = 1/2$ is of particular interest because it is analogous to the quadratic split \eqref{eq:quadSplit} used to create stable split forms for the compressible Euler equations. This form is considered here for the Burgers' equation in an attempt to anticipate the behaviour of relevant split form DGSEM as applied to Euler-based turbulence, which will be addressed in the remaining sections of the present work through the inviscid Taylor-Green vortex problem. The conservative and non-conservative forms of \eqref{eq:split_quads} are also considered due to their obvious significance. These three forms are compared in Fig.~\ref{fig:1dComparison}, which shows the transformed (modal space) values of the right-hand side of \eqref{eq:split_quads}, \comment{namely $\mbox{TRHS}_i \approx \partial_t \hat{q}_i = \mathbf{V}^{-1} \partial_t \check{q}$, for different modal coefficients $i = 0, \dots, N$.}

We stress that TRHS values do not take into account boundary contributions, as they are intended to represent the effect of the volumetric terms on the numerical solution. \comment{Also, for the one-dimensional test addressed, the interface flux contribution does not depend on $Q$ or even $\alpha$.} The expansion order, $N$, varies for the different plots in Fig.~\ref{fig:1dComparison}. The exact values of $\mbox{TRHS}_i$ have been computed via consistent integration of the conservative term in \eqref{eq:split_quads} and are also shown in the plots for reference. We note that even higher values of $N$ were considered and the main trends remained the same. Random perturbations of up to 50\% in the coefficients of \eqref{eq:modal_expansion} have been added, and again the main tendencies observed in Fig.~\ref{fig:1dComparison} remained. This demonstrates that the results obtained are not caused by a fortuitous, particular choice of parameters. \comment{Note also that the overall profile shown in Fig.~\ref{fig:1dComparison} is physically consistent in the sense that energy is flowing from low-order modes (low frequencies) to the high-order ones (high frequencies), as happens in the energy cascade mechanism typical of turbulence.}

\comment{The most relevant observation one can draw from Fig.~\ref{fig:1dComparison} is that conservative DGSEM tends to over-predict (intensify) the coefficients' variation rate $\mbox{TRHS}_i \approx \partial_t \hat{q}_i$ whereas non-conservative DGSEM under-predicts it, this being especially significant at the highest-order modes. In terms of modal energy, one can expect these results to be even more significant for the relative variation rates $(\partial_t \hat{q}_i^2) / \hat{q}_i^2$, which, in the frozen solution analysis considered, should scale like
\begin{equation}
	\frac{1}{\hat{q}_i^2} \frac{\partial \hat{q}_i^2}{\partial t} \propto \frac{1}{\hat{q}_i^2} \left( \hat{q}_i \frac{\partial \hat{q}_i}{\partial t} \right) \approx \frac{\mbox{TRHS}_i}{\hat{q}_i} \propto (i+1)^{5/6} \, \mbox{TRHS}_i \mbox{ .}
\end{equation}
In other words, higher-order modes suffer a larger relative deviation (over/under-prediction) in terms of modal energy variation rate since they are inherently less energetic. As a result, standard conservative and non-conservative DGSEM approaches without over-integration are expected to grossly miscalculate the dynamics of small turbulent scales, eventually inducing numerical instabilities.}

The over-prediction of modal energy at higher modes in conservative DGSEM approaches has been reported in previous studies \cite{Kirby2003} and correctly attributed to polynomial aliasing errors due to inexact integration of the non-linear terms. Here we have shown that the quadratic split form tends to predict energy variations more accurately due to the averaging of opposite conservative and non-conservative tendencies. In a sense, this can be regarded as a built-in dealiasing mechanism as it suppresses aliasing errors. This result is similar to that reported for pseudo-spectral methods \cite{blaisdell1996effect}, although the reasons are possibly different because aliasing errors here originate from inexact polynomial quadratures, whereas in \cite{blaisdell1996effect} the aliasing errors stem from the handling of Fourier transforms with an inconsistent number of modes.

Another relevant point is that energy variation levels yielded here by split form DG approaches tend to be slightly below the correct levels, which can be advantageous for robustness and should help to keep energy ``under control," suppressing excessive growth of small-scale energy and preventing numerical divergence in under-resolved turbulence simulations. It has also been found (not shown) that when the number of quadrature nodes $Q$ is increased for the conservative form, energy variation levels approach the exact values monotonically from above, as over-integration gradually reduces the TRHS levels over-predicted by the conservative form and bring them to the correct values. This qualitatively explains why split forms can be more stable than conservative discretisations even when quadratures are consistent (with non-linearities of the problem) albeit still inexact. \comment{This is in line with a recent study \cite{manzanero2017insights} that proved over-integration unable to entirely eliminate aliasing errors in advection problems with non-constant speed functions, whereas appropriate split forms were able to.}

\comment{We recall that when the compressible Euler equations are written using conserved variables, the flux functions to be integrated numerically involve rational functions which, in general, cannot be integrated exactly with a finite number $Q$ of quadrature nodes. Integrations are however expected to converge to the exact values as $Q$ increases, which could grant higher accuracy for over-integration approaches at high enough $Q$, if they are stable (and affordable). The accuracy and robustness of split and over-integrated conservative forms as applied to under-resolved turbulence are assessed and compared in the following sections.}

\section{TGV simulations and numerical stability}
\label{sec:analysisTools}

The Taylor-Green vortex flow was introduced in \cite{taylor1937mechanism} as a model problem for the analysis of transition and turbulence decay. The test case was originally proposed for the incompressible Navier-Stokes equations in a cubic domain with (triply-)periodic boundary conditions. As in previous studies \cite{shu2005numerical,drikakis2007simulation}, we adopt a modified version of the initial conditions, which is suited for compressible flow solvers. The following expressions have been used as the initial state within $\Omega = [-\pi \ell_o, \pi \ell_o]^3$, respectively, for the density, the three velocity components, and the static pressure:
\begin{gather}
	\rho = \rho_o \mbox{,} \label{eq:TGV1}\\
	u = V_o \sin{(x/\ell_o)} \cos{(y/\ell_o)} \cos{(z/\ell_o)} \mbox{,} \,\,\,\,
	v = - V_o \cos{(x/\ell_o)} \sin{(y/\ell_o)} \cos{(z/\ell_o)} \mbox{,} \,\,\,\,
	w = 0 \mbox{,} \\
	p = \rho_o c_o^2 / \gamma + \rho_o V_o^2 \left[\cos{(2x/\ell_o)}+\cos{(2y/\ell_o)}\right] \left[2+\cos{(2z/\ell_o)}\right] / 16 \mbox{,}\label{eq:TGV5}
\end{gather}
the total energy per unit volume being $E = p/(\gamma-1)+\rho(u^2 + v^2 + w^2)/2$. For the reference quantities, we have adopted the values $\ell_o = \rho_o = V_o = 1$ and $c_o = 10$, leading to a Mach number of $0.1$, which makes our cases nearly incompressible. A non-dimensional time $t$ is adopted based on the scale $\ell_o / V_o = 1$. A Reynolds number could be defined as $Re = \rho_o V_o \ell_o / \mu_o$, but only the inviscid problem is considered in this study, where the compressible Euler equations have been simulated directly, with $\gamma = 1.4$. \comment{Similar to Gassner et al. \cite{gassner2016split} we take the final time of the TGV flow evolution to be $T=14$.}

The Euler-based simulations considered in the present study are expected to be representative of viscous TGV solutions at very high Reynolds numbers, provided that numerical dissipation mimics the dissipative character of Navier-Stokes turbulence in the limit of vanishing viscosity \cite{frisch1995turbulence}. The equivalent scenario in traditional LES would be to have zero molecular viscosity and rely solely on the regularization of a subgrid-scale model. Here, upwind dissipation is expected to play this latter role. In \cite{moura2016ontheeddy}, the inviscid TGV was simulated with different Riemann solvers in the context of DGSEM with consistent integration. The study showed in particular that Roe-based discretisations were more robust than LLF-based ones, which crashed more easily - a surprising result in itself, which will be discussed below. However, at higher polynomial orders ($N \geq 6$) both fluxes produced test cases that lacked stability. Subsequent DGSEM simulations \cite{gassner2016split} based on suitable split form discretisations demonstrated remarkable robustness and were capable of yielding stable solutions of the inviscid TGV with both Roe and LLF-type fluxes, even at very high-orders (e.g.\ $N = 15$). In the following we summarize the test cases addressed in \cite{moura2016ontheeddy,gassner2016split}, consider new (lower-order) test cases, and discuss the inherent instability issue of the inviscid TGV flow.

The base set of test cases considered are given in Table\,\ref{table:rcm1}. Each column corresponds to the number of one-dimensional polynomial modes $m = N + 1$ used, $N$ being the polynomial order. The number of elements varies by a factor of $\sqrt{2}$ between adjacent rows to reduce computational costs, as choosing a factor of 2 would requires a large number of degrees of freedom in the last line. The values of $n_{el}$ were chosen so that the degrees of freedom $N_{dof} = (n_{el} \, m)^3$ are kept (approximately) constant along a given row. Equispaced grids of cubic elements have been employed. Values in the Table's core represent the number of elements in each spatial direction, $n_{el}$. Crossed out numbers denote simulations that crashed with consistent integration. The KG and DU split form discretisations were however able to stabilize all the tests cases, see \cite{gassner2016split} for a complete discussion on the robustness of other available split forms. Spectral element codes $Nektar$++ \cite{cantwell2015nektar} and $FLUXO$ (\texttt{http://www.github.com/project-fluxo}) have been used respectively for the computations based on consistent integration and split forms.

\begin{table}[h]
	\caption{Summary of cases --- crossed out cases crashed with consistent/over- integration, whereas all test cases ran to the final time with the KG and DU split forms.}
	\centering
	\begin{tabular}{c ccccccc c ccccccc}
		\\ [-2.0ex]
		\hline \\ [-2.0ex]
		\text{ } &
		\multicolumn{7}{c}{Roe} &
		$\!\!\!\!$ &
		\multicolumn{7}{c}{LLF} \\
		$m = N + 1$ &
		2 & 3 & 4 & 5 & 6 & 7 & 8 &
		$\!\!\!\!$ &
		2 & 3 & 4 & 5 & 6 & 7 & 8 \\
		\hline \\ [-2.0ex]
		\text{ } &
		56 & 37 & 28 & 23 & 19 & 16 & 14 &
		$\!\!\!\!$ &
		56 & 37 & 28 & 23 & \xcancel{19} & \xcancel{16} & \xcancel{14} \\
		$n_{el}$ &
		79 & 52 & 39 & 32 & 28 & 23 & \xcancel{19} &
		$\!\!\!\!$ &
		79 & 52 & 39 & 32 & \xcancel{28} & \xcancel{23} & \xcancel{19} \\
		\text{ } &
		112 & 75 & 56 & 45 & 39 & \xcancel{32} & \xcancel{28} &
		$\!\!\!\!$ &
		112 & 75 & 56 & 45 & \xcancel{39} & \xcancel{32} & \xcancel{28} \\
		\hline
	\end{tabular}
	\label{table:rcm1}
\end{table}

As pointed out in \cite{moura2016ontheeddy}, the crashes are probably not related to insufficient integration. We acknowledge that exact integration is never achieved with a finite number of quadrature points since the flux vector components are not polynomials, but rational functions, given that our consistent integration approach is based on the conservative form of the Euler equations. However, since quadrature errors still tend to zero as the number of integration points is increased \cite{trefethen2013approximation}, they are expected here to be negligible in the case of consistent integration, given the low Mach number adopted. In Table\,\ref{table:rcm2}, we consider the time of crash, $t_c$, versus the number of integrations points $Q$ employed (per element and per dimension) for test case $m = 8$, $n_{el} = 14$. Consistent integration is performed both for the volume and surface terms, which take $Q^3$ and $Q^2$ points, respectively. We note that this scaling is what makes consistent integration prohibitive at high polynomial orders. The asymptotic behaviour of $t_c$ for LLF (particularly as $Q$ is increased from 16 to 32) indicates that consistent/over- integration has already contributed all it could towards stabilisation. \comment{Further, we note that increasing the value of $Q$ does offer stabilisation for Roe and entries in Table \ref{table:rcm2} marked with ``---'' ran successfully.}

\begin{table}[h]
	\caption{Time of crash $(t_c)$ vs.\ quadrature points $Q$ for test case $m = 8$, $n_{el} = 14$}
	\centering
	\begin{tabular}{c ccccccccc}
		\\ [-2.0ex]
		\hline \\ [-2.0ex]
		$Q$ & 9 & 10 & 11 & 12 & 13 & 14 & 15 & 16 & 32 \\
		\hline \\ [-2.0ex]
		LLF $t_c$ & 4.5434 & 6.4498 & 6.6902 & 8.4094 & 8.3950 & 8.3930 & 8.3932 & 8.3953 & 8.3954 \\
		Roe $t_c$ & 4.4824 & 7.4659 & 9.0643 & --- & --- & --- & --- & --- & --- \\
		\hline
	\end{tabular}
	\label{table:rcm2}
\end{table}

For Legendre-Gauss-Lobatto quadratures, the number of nodes required for the consistent integration of linear, quadratic and cubic terms are respectively $Q \approx m$, $Q \approx 3m/2$ and $Q \approx 2m$. 
We see from Table \ref{table:rcm2} that LLF's $t_c$ increases with $Q$ until up to $Q = 3m/2$ and then remains practically unaffected. The same is true for the Roe flux, except that Roe-based cases are stabilized for $Q \geq 3m/2$. This is consistent with the nearly incompressible nature of the TGV cases considered, since density variations are small and the terms being integrated are essentially quadratic. Despite of this, consistent integration of all the cases in Table\,\ref{table:rcm1} assumed a cubic non-linearity for the compressible Euler equations and used $Q = 2m$. The differences between LLF and Roe are likely due to strong over-upwinding effects induced by the former, which increase the likelihood of TGV instabilities, as explained in the following.

There is an open debate in the literature as to whether or not the inviscid TGV flow might develop singularities which lead to the actual collapse of the solution \cite{brachet1992numerical,cichowlas2005evolution,hou2008blowup,gibbon2008three}. However, as emphasized in \cite{moura2016ontheeddy}, this possibility is only considered for the exact, energy conserving solution of the Euler equations. This is in contrast with the character of Navier-Stokes turbulence in the limit of vanishing viscosity, which remains dissipative and is not expected to develop singularities, see \cite{frisch1995turbulence} Secs. 5.2 and 9.3. We stress that any LES-like approach should follow the latter behaviour when viscosity is set to zero as subgrid dissipation (explicit or implicit) remains active. Nevertheless, an energy-conserving bias has been found to be partially induced by (some of) the discretisations considered due to the very sharp spectral dissipation expected at higher orders \cite{moura2015} and also from over-upwinding \cite{moura2016anlesICOSAHOM}. It is known that a sharp dissipation in Fourier space can induce an energy-conserving behaviour in turbulent simulations \cite{banerjee2014transition}. This is because, under certain circumstances, a cutoff-like spectral dissipation may prevent the formation and subsequent destruction of small scales near the cutoff wavenumber, acting effectively as a barrier to the inertial cascade and causing a pile-up of small-scale energy. The LLF solver is expected to have a particularly sharp dissipation due to its over-upwind bias for the momentum equations, owing to the disparity between acoustic and convective wave speeds, especially at low Mach numbers, see \cite{moura2016anlesICOSAHOM}. As explained in \cite{moura2016ontheeddy}, the main evidence that this behaviour is taking place (even if partially) is the so-called ``energy bump" observed in the energy spectra of LLF-based DGSEM simulations, see e.g., \cite{diosady2015higher,wiart2015implicit}. It has been demonstrated \cite{cichowlas2005effective,frisch2008hyperviscosity,banerjee2014transition} that pre-dissipative bumps emerge as the solution begins to follow an energy-conserving dynamics when only a finite number of Fourier modes are retained (limit of increasingly sharp dissipation). More evidence of this artificial energy-conserving character will be given in the next section.

Besides the occurrence of the energy-conserving bias mentioned above, it's been suggested in \cite{moura2016ontheeddy} that simulations following this character (higher-order discretisations, especially LLF-based ones) are also inducing the emergence of the physical singularities long conjectured for the inviscid energy-conserving TGV problem. This second point is however somewhat speculative and may be very difficult to prove. Still, note that while suitable split form discretisations are able to prevent instabilities with both Roe and LLF solvers even at very high orders, fully central split forms which are nearly kinetic energy conserving are only stable at low polynomial orders $(N \leq 3)$. This further supports our claim that the TGV instabilities observed are likely of physical origin, as their likelihood of emergence increases in higher-order (better resolved) discretisations of energy-conserving character.

In summary, although the instabilities observed are not entirely understood at this point, it is clear that suitable split form discretisations have superior non-linear stability characteristics. Robustness, however, can be easily obtained at the cost of solution quality, and so the accuracy of split form DG approaches require further investigation. Note that such solutions are still prone to aliasing errors because they do not rely on standard consistent/over- integration. The evaluation of the accuracy and fidelity of suitable split form DG approaches by comparison to consistent/over- integration for the invisicid TGV flow is the main focus of the next section. To the authors' knowledge this is the first comparison between these two dealiasing approaches to be conducted in the context of under-resolved turbulence computations.

\section{Solution quality for different dealiased DGSEM}
\label{sec:solQuality}

This section is devoted to the assessment of solution quality of the test cases introduced in the previous section.  
We investigate different polynomial orders and number of DOFs from the configurations in Table \ref{table:rcm1} originally used to examine robustness. To assess and compare the solution quality of consistent/over- integration and split forms we select the KG split as a representative scheme. This is motivated by previous results applying split form DGSEM to the inviscid TGV vortex. Specifically, while the most stable schemes considered in \cite{gassner2016split} yielded similar flow solutions, the KG split combined a simple implementation with reduced computational cost. However, for completeness, we provide a brief comparison of the KG and DU split forms in \ref{sec:compSplit} to further justify the choice of KG as a representative split form method for the compressible Euler equations.

We first examine the effect of the Riemann solver at moderately high polynomial orders in Section~\ref{sec:effectRiemann}, both for consistent/over- integration and split forms. Then, Section~\ref{sec:strangeQuality} discusses some results obtained at very high-orders and also the behaviour of low-order discretisations and provides direct comparisons between consistent/over- integration and split form discretisations.

\subsection{Effect of the Riemann flux at moderately high-orders}
\label{sec:effectRiemann}

\begin{figure}[!ht]
	\centering
	\subfloat[$m = 5$, $n_{el} = 32$: Roe numerical flux.] 
	{
		\includegraphics[scale=0.60]{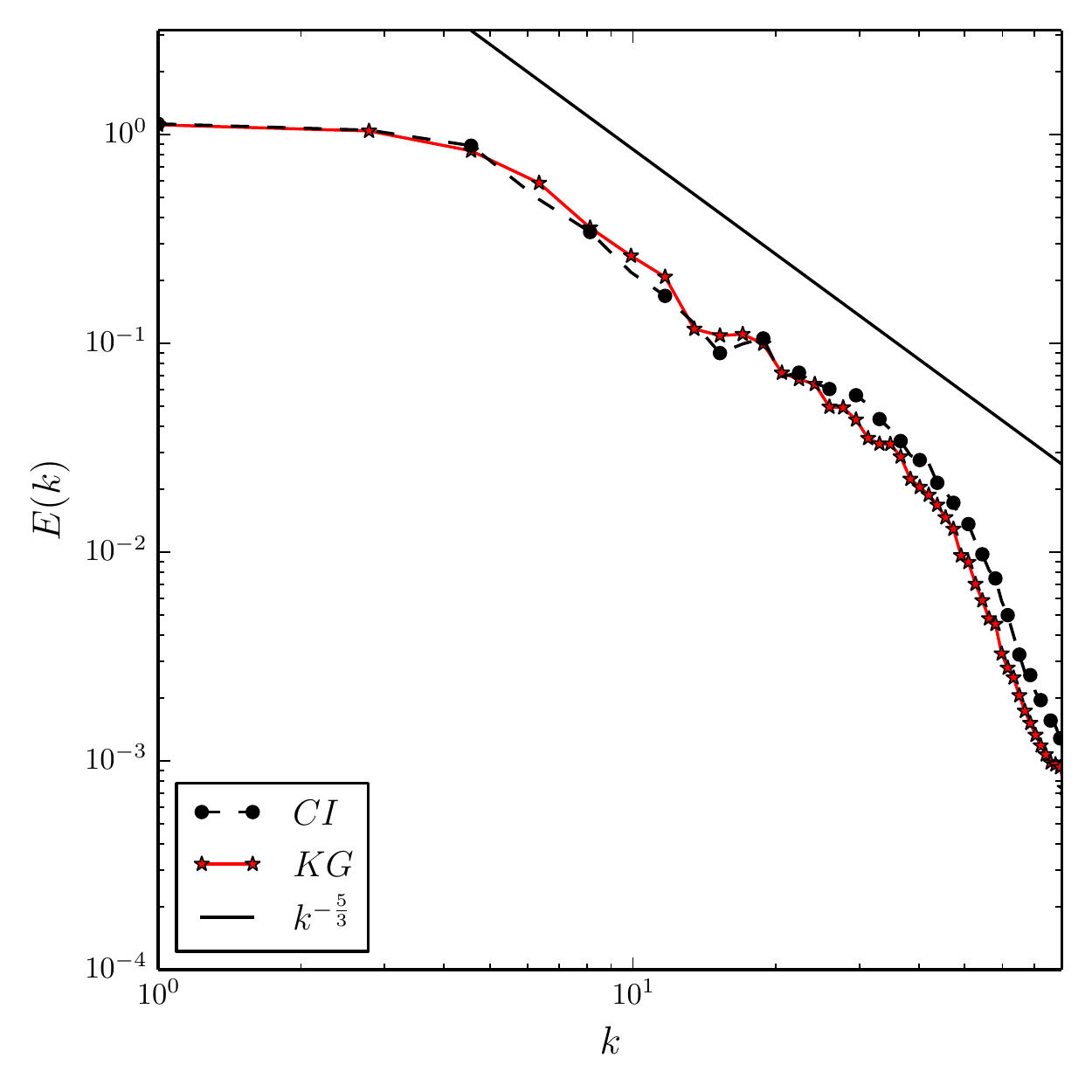}
	}
	\subfloat[$m = 5$, $n_{el} = 32$: LLF numerical flux.] 
	{
		\includegraphics[scale=0.60]{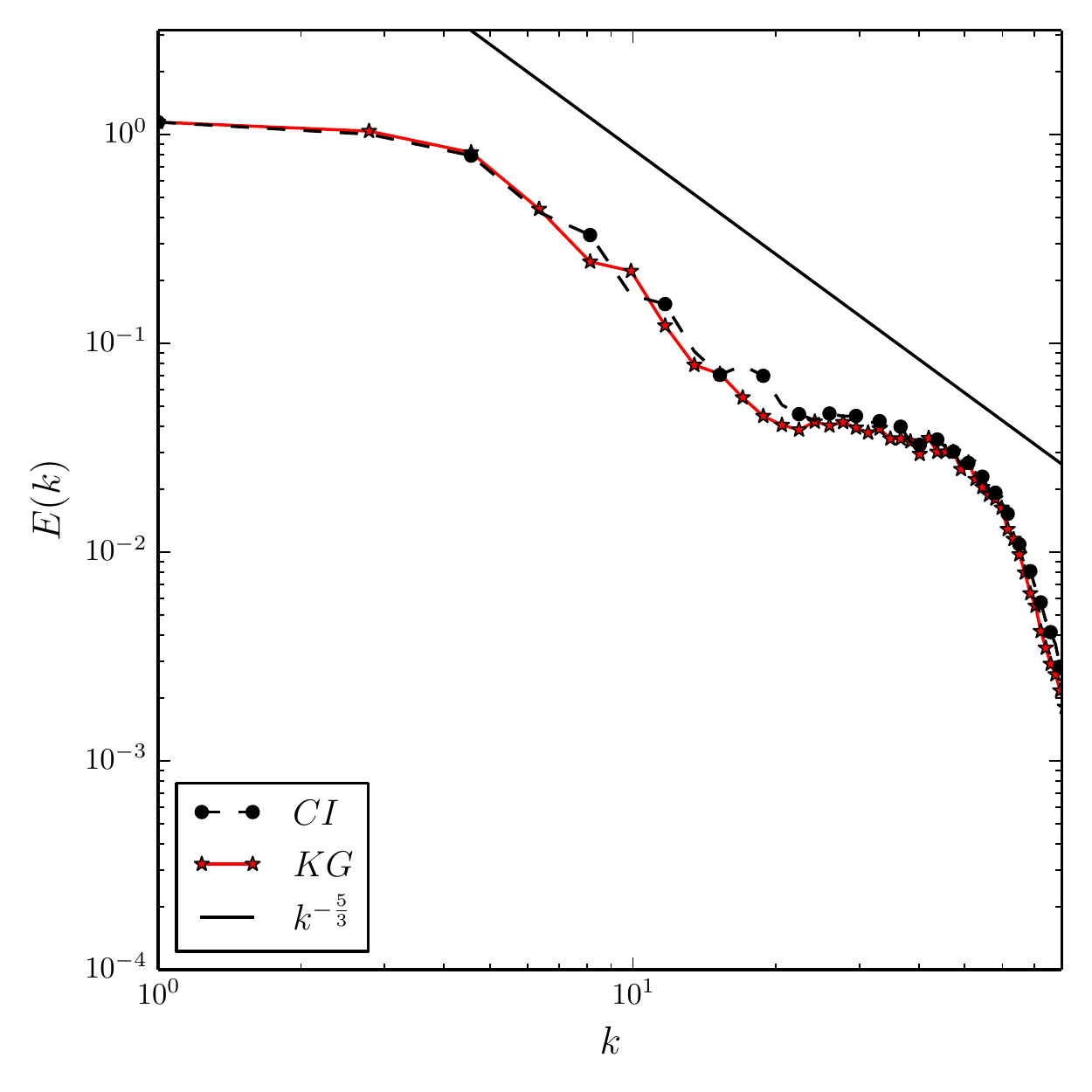}
	}
	\caption{Comparison between energy spectra of consistent integration (CI) and split form (KG) for case $m = 5$, $n_{el} = 32$. Comparing the two plots also makes evident the spurious energy bump in LLF-based cases.} 
	\label{fig:X2}
\end{figure}

High-order results ($m \geq 4$) obtained with consistent integration have been partially analysed in \cite{moura2016ontheeddy}. A complete accuracy assessment is not possible as a DNS solution for the inviscid TGV is simply out of reach. One of the main conclusions in \cite{moura2016ontheeddy} was that LLF-based discretisations produced solutions with an excess of small-scale energy owing to the formation of an ``energy bump," as mentioned in Section~\ref{sec:analysisTools}. Roe-based solutions, on the other hand, did not show this spurious feature and produced energy spectra similar to those from classic LES of isotropic turbulence at infinite Reynolds number (i.e.\ with molecular viscosity set to zero). Here, these results have also been observed for the split forms considered. More specifically, KG energy spectra have been found to be very close to those obtained with consistent integration. This is shown in Fig.~\ref{fig:X2} for the test case $m = 5$, $n_{el} = 32$.

It is interesting to note from Fig.~\ref{fig:X2} that, at small scales (say, beyond the centre of the inertial range), split form spectra are slightly less energetic than those obtained with consistent integration. This has been observed as a general trend in our comparisons and is consistent with the estimates of Section~\ref{sec:aliasingProperties}. The close proximity between results based on split forms and consistent integration have been found to hold also at higher orders (cf.\ Fig.~\ref{fig:X4}(a)). However, for low-order discretisations ($m \leq 3$), results exhibited non-negligible differences, which is also consistent with Section~\ref{sec:aliasingProperties}. Low-order solutions will be discussed in Section~\ref{sec:strangeQuality}, along with some results obtained at very high-orders.

\begin{figure}[!ht]
	\centering
	\subfloat[$m = 5$, $n_{el} = 32$: CI LLF numerical flux.]
	{
		\includegraphics[scale=0.505]{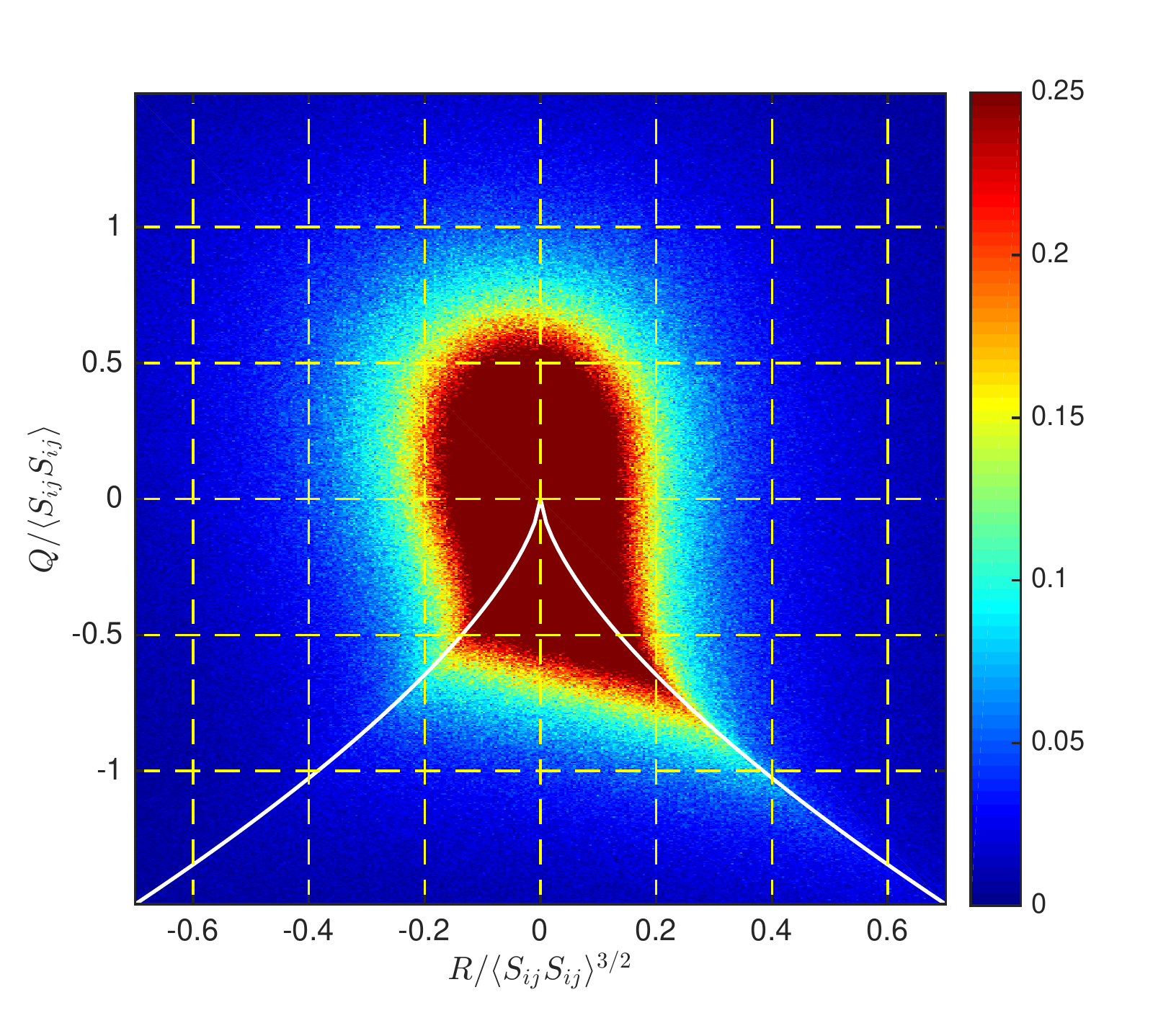}
	}
	\subfloat[$m = 5$, $n_{el} = 32$: CI Roe numerical flux.]
	{
		\includegraphics[scale=0.505]{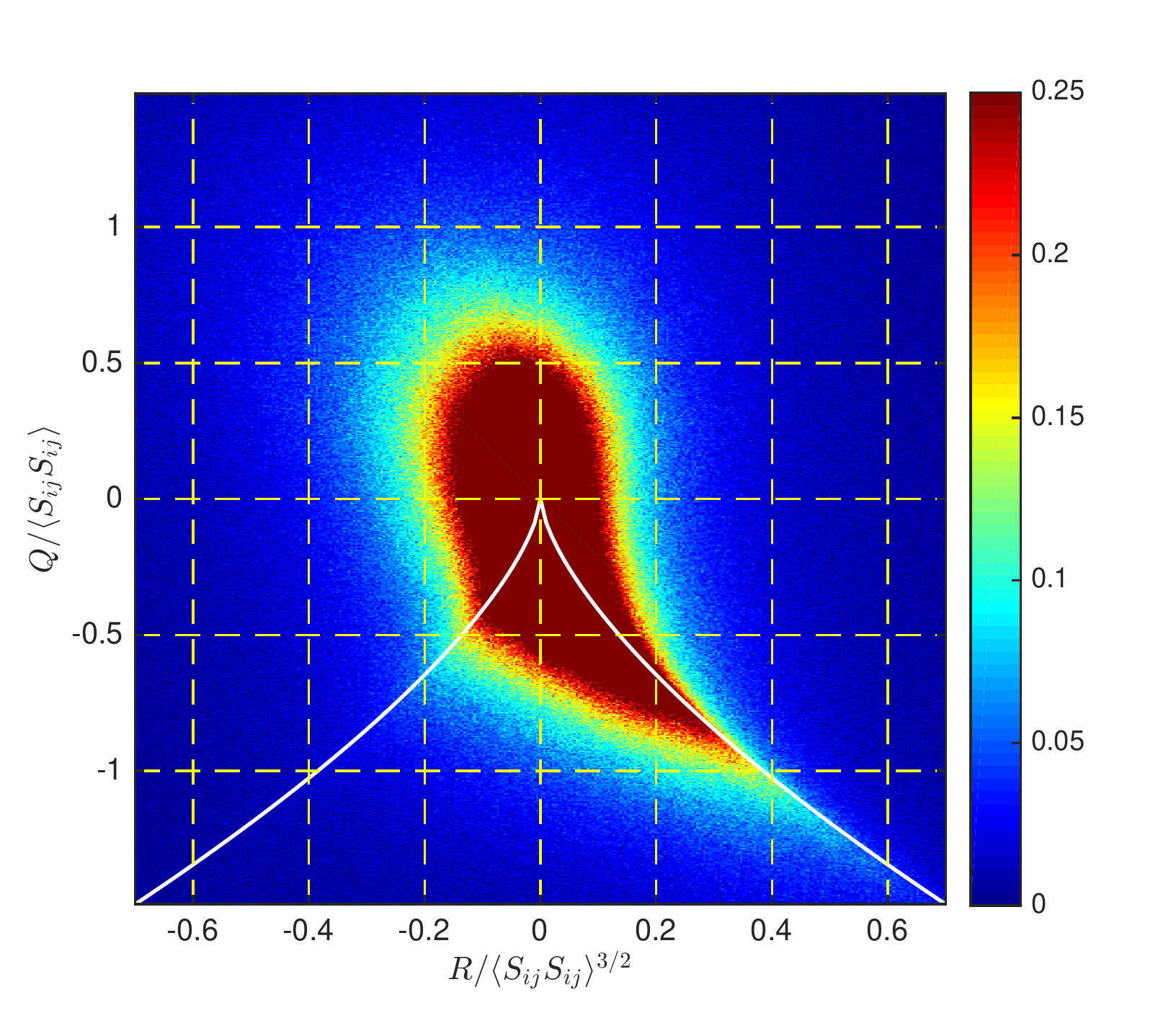}
	}
	\\
	\subfloat[$m = 5$, $n_{el} = 32$: KG LLF numerical flux.]
	{
		\includegraphics[scale=0.505]{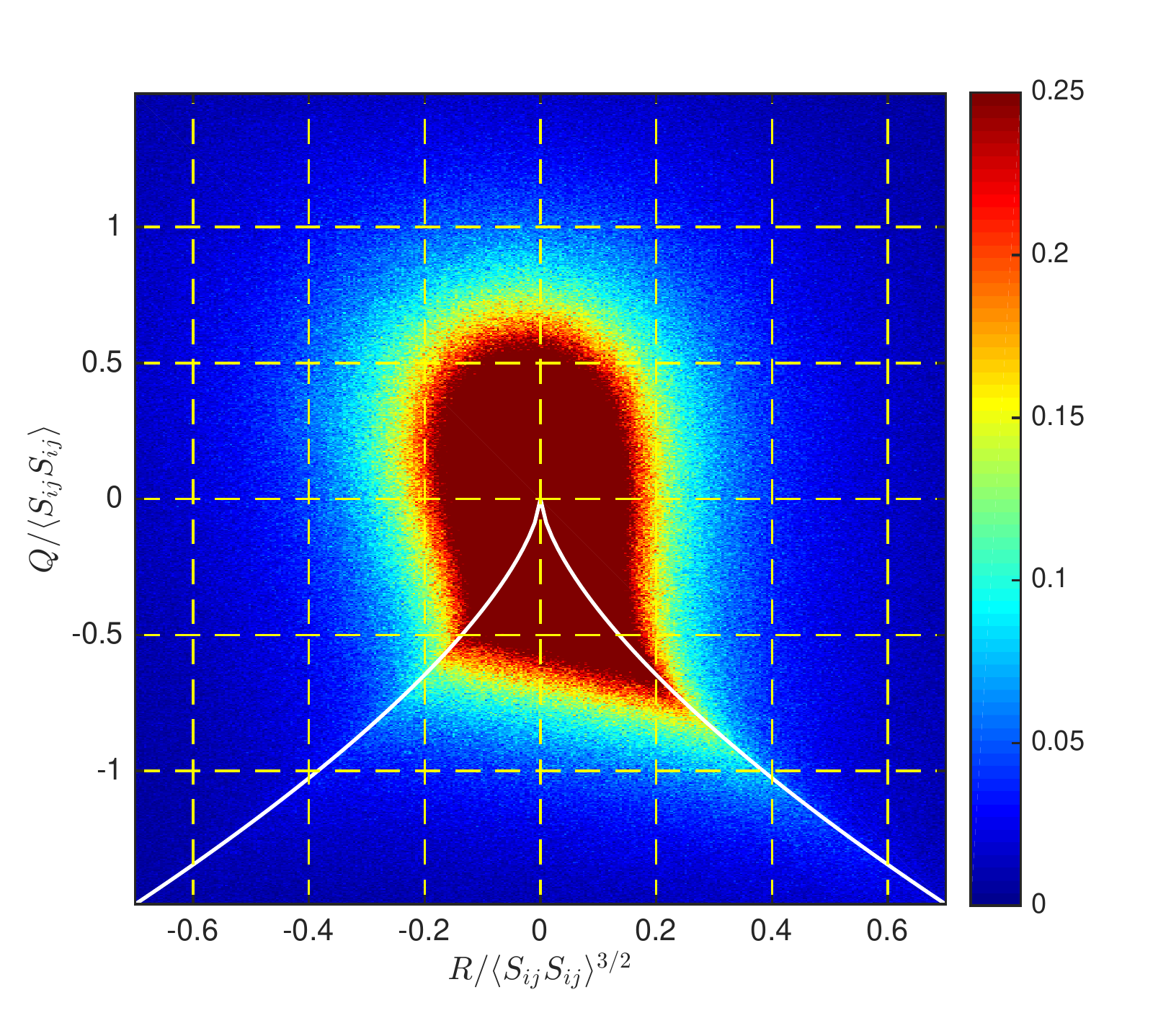}
	}
	\subfloat[$m = 5$, $n_{el} = 32$: KG Roe numerical flux.]
	{
		\includegraphics[scale=0.505]{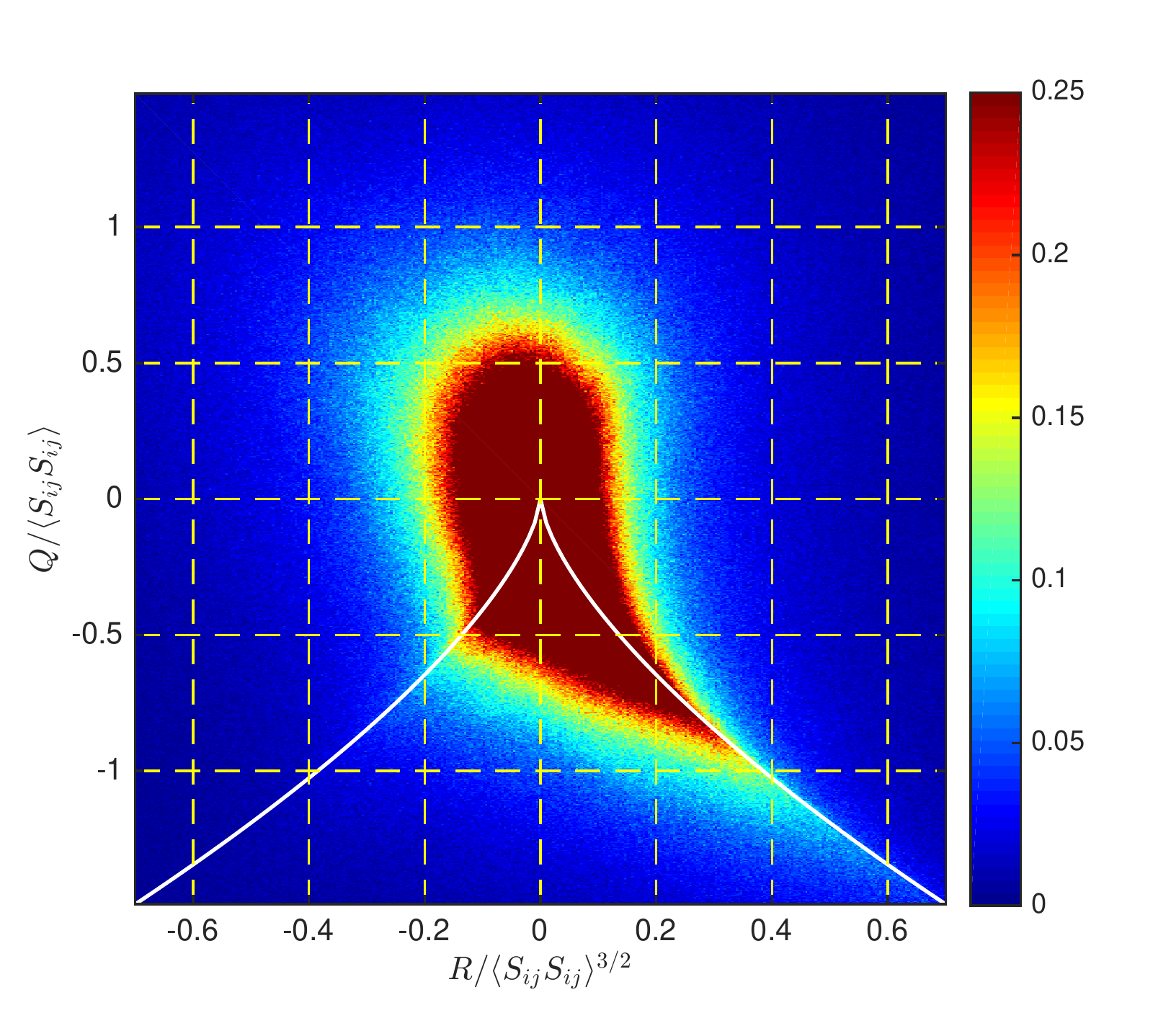}
	}
	\caption{QR diagrams of case $m = 5$, $n_{el} = 32$ at $t = 9$ for consistent integration (top) and split forms (bottom) with LLF (left) and Roe-type (right) numerical fluxes.}
	\label{fig:X3}
\end{figure}

The results discussed so far indicate that high-order split form discretisations can provide solutions of quality very similar to that of consistent integration. An additional assessment of solution quality for both discretisations is presented in Fig.~\ref{fig:X3}, namely, the QR diagrams of the solutions whose spectra have been considered in Fig.~\ref{fig:X2}. These consist of joint PDFs of invariants Q and R of the velocity gradient tensor, which is evaluated at a large set of points within the flow domain. The resulting diagram provides specific insight about the kinematic topology of the flow \cite{Chong1990}. More details about QR diagrams and their construction is given in Appendix B. Typical QR diagrams exhibit a characteristic teardrop profile, as seen, e.g., in Fig.~\ref{fig:X3}(b), in various turbulent scenarios, which is considered a qualitatively universal feature of turbulence \cite{tsinober2009aninformal}. In particular, a prolonged extension of the profile over the bottom-right quadrant of the diagram, known as the ``Viellefosse tail'' \cite{Viellefosse1984}, is expected. The canonical shape for a teardrop profile of isotropic turbulence can be found in e.g. \cite{Laizet2015}. More generally, the overall profile shape observed in a QR diagram can reveal how closely numerical results are able to reproduce physical turbulent features. Therefore, it is a good indicator to assess the quality of the simulations run for the two dealiasing strategies under consideration. It is also useful to check how the lack of resolution or the polynomial order of the approximation can affect the QR diagram.

An interesting result is shown in Fig.~\ref{fig:X3}, namely, that Roe-based QR profiles very closely resemble the canonical teardrop shape \cite{Laizet2015}, whereas LLF-based ones do not conform so well to that shape and feature a much less pronounced tail. This further supports that LLF-based solutions are less physical than Roe-based ones and indicates that energy bumps contribute to the overall turbulent dynamics with kinematic states that are artificial. In fact, because LLF-based QR profiles are more symmetrically oriented with regards to the origin of the diagrams, bump-related kinematic states are probably somewhat random, while LLF-based QR profiles even resemble an artificially generated (Gaussian) turbulent state considered in \cite{tsinober2009aninformal}.

Assuming that energy bumps are in fact emerging from an energy-conserving dynamics (as discussed in Section 4), theory predicts that bump-related scales should reach an equilibrium (loosely called a ``thermalised" state) where energy equipartition is favoured \cite{frisch2008hyperviscosity}, which is consistent with the more symmetrical distribution of kinematic states of LLF-based QR profiles. This hypothesis will receive further confirmation in Section~\ref{sec:strangeQuality}, where very high-orders (with sharper dissipative behaviours) are considered. A final remark from Fig.~\ref{fig:X3} is that QR profiles of split forms are only slightly wider than those of consistent integration, but otherwise practically indistinguishable. This supports our claim that high-order split form discretisations provide solutions of quality very similar to that of consistent/over- integration for the inviscid TGV flow.

\subsection{Solution quality at lower and higher polynomial orders}
\label{sec:strangeQuality}

\begin{figure}[!ht]
    \centering
	\subfloat[$m = 8$, $n_{el} = 14$: Roe numerical flux.]
    {
        \includegraphics[scale=0.6]{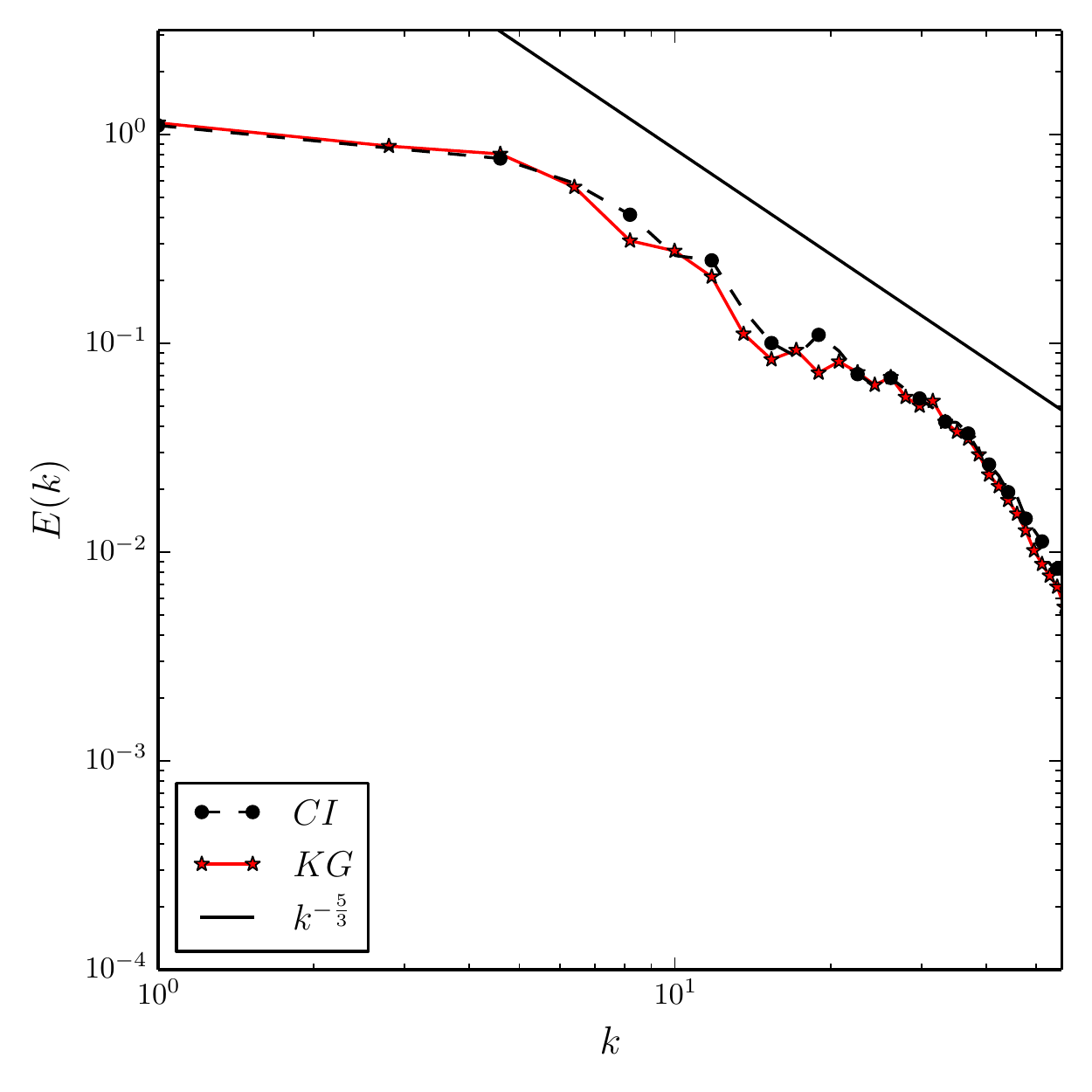}
    }
    \subfloat[$m = 16$, $n_{el} = 7$: Roe numerical flux.]
    {
        \includegraphics[scale=0.6]{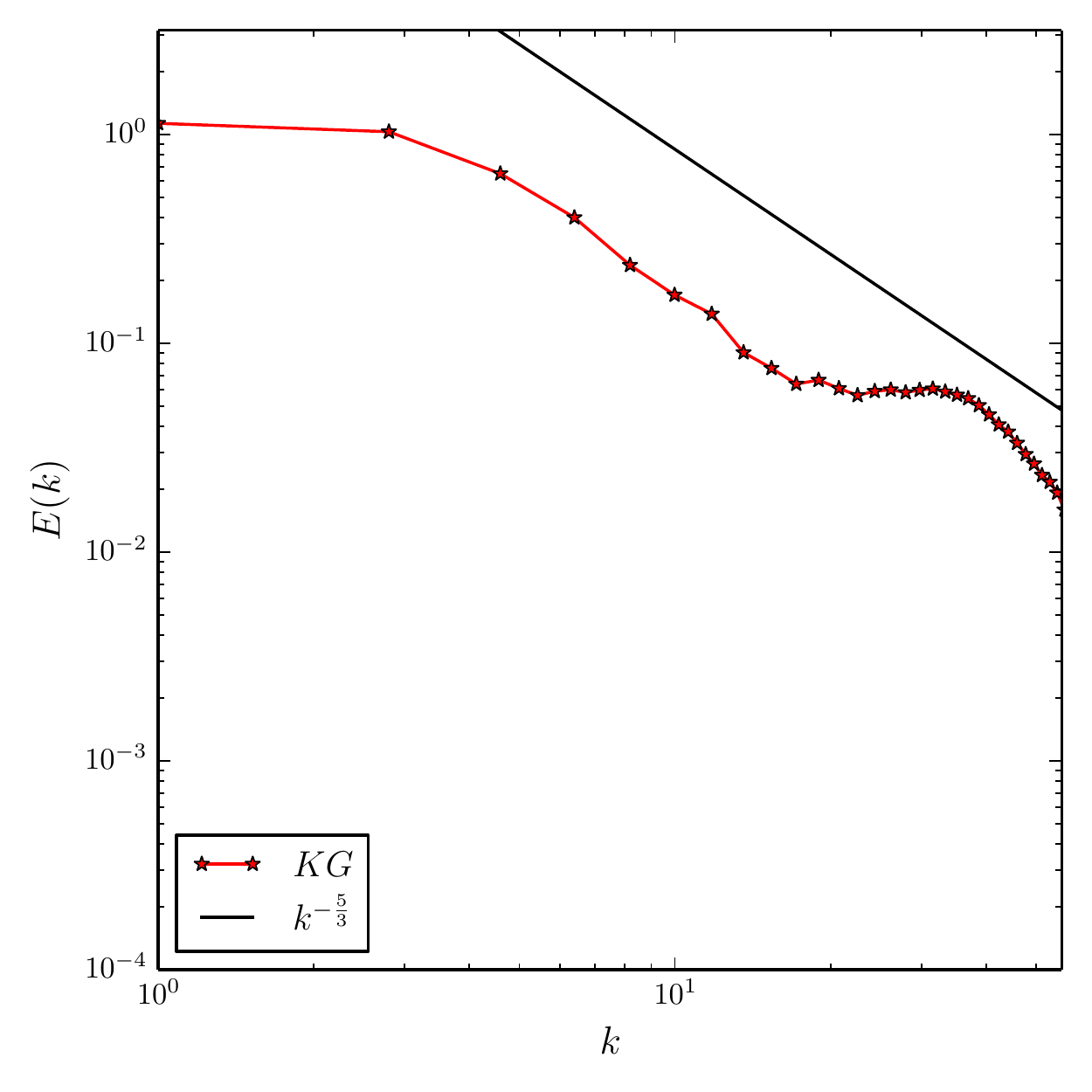}
    }
    \caption{Energy spectra at $t = 9$ obtained with Roe-type numerical flux for same-DOF cases with (a) $m = 8$ and (b) $m = 16$.}
    \label{fig:X4}
\end{figure}

We start by considering Roe-based solutions at higher orders, namely $m=8$ and $m=16$ in Fig~\ref{fig:X4}. Note that these have the same number of DOFs and that the spectrum of the equivalent test case at $m=5$ can be seen in Appendix A, cf.\ Fig.~\ref{fig:X1}(b). At $m=16$, consistent integration is not sufficient to suppress TGV instabilities and hence only the result from the KG scheme is shown here. One can see that an energy bump can emerge even with Roe-type numerical flux as the polynomial order is increased. This is consistent with our claim that energy bumps are caused by a sharper dissipation in wavenumber space. Also, it seems from Fig.~\ref{fig:X4} that the higher order solutions only achievable with split forms follow the trends that would have been obtained with consistent/over- integration if these simulations were stable.

\begin{figure}[!ht]
    \centering
    \subfloat[$m = 16$, $n_{el} = 7$: Roe numerical flux.]
    {
        \includegraphics[scale=0.505]{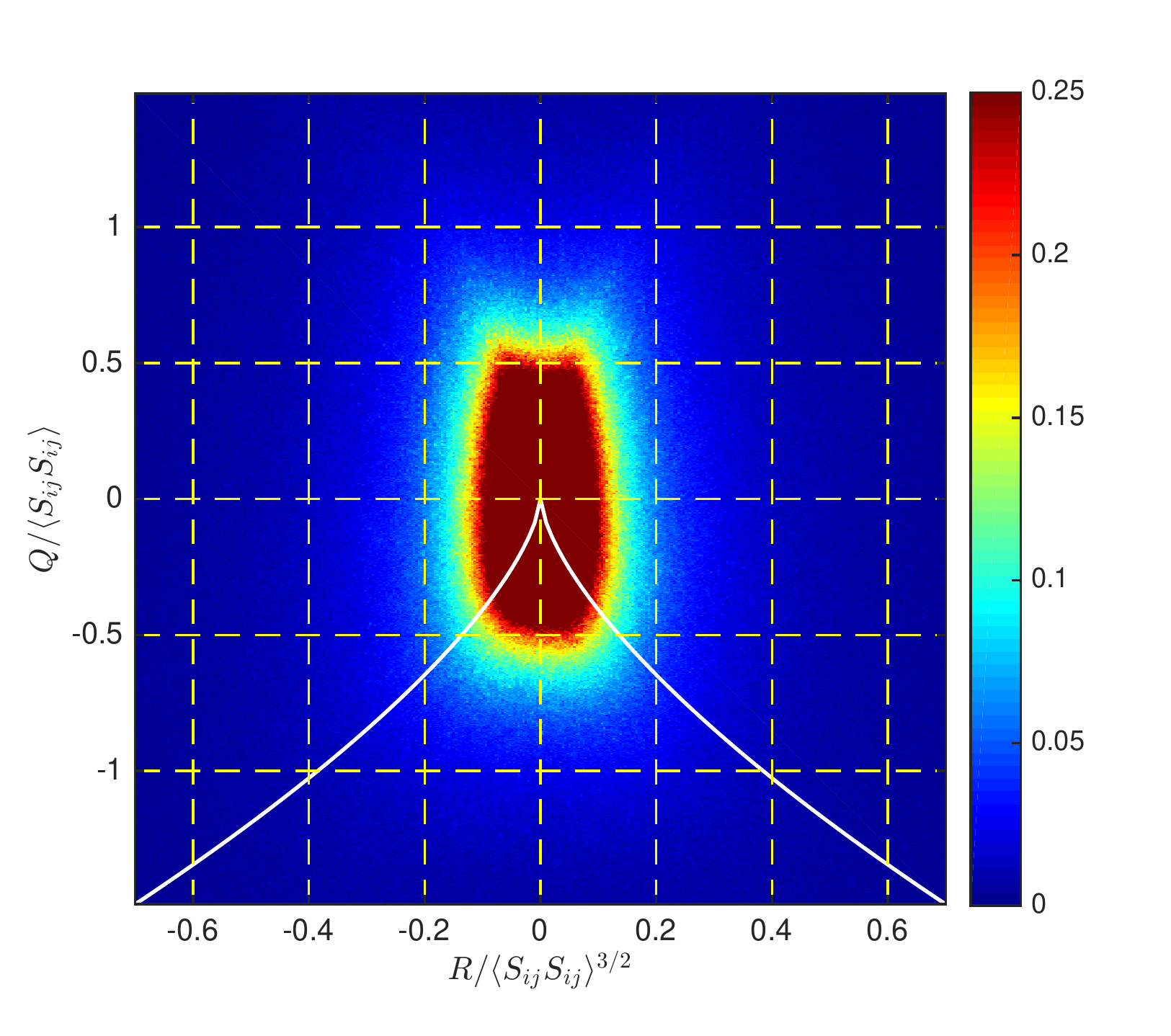}
    }
    \subfloat[$m = 4$, $n_{el} = 28$: Fully central KG form.]
    {
        \includegraphics[scale=0.505]{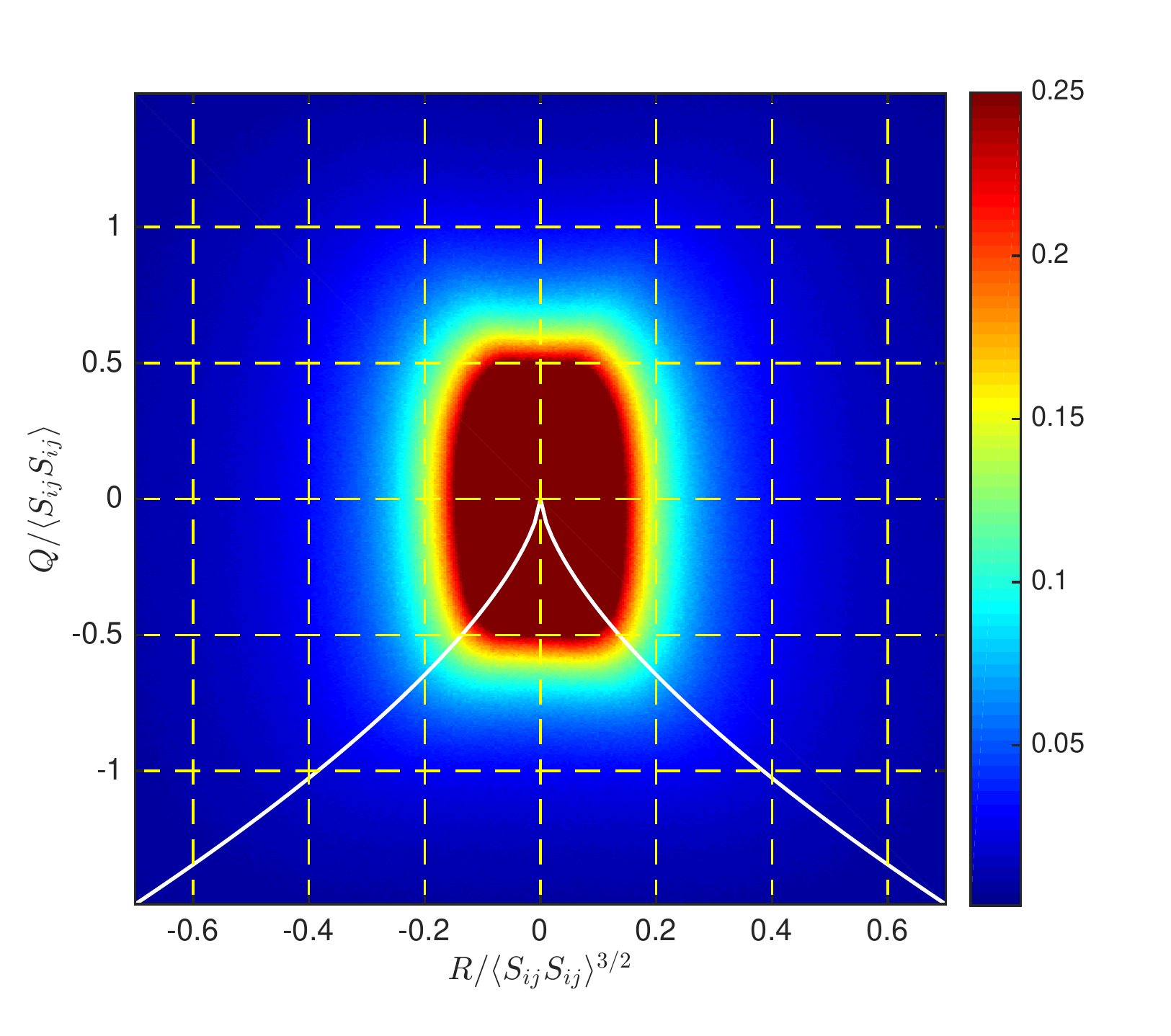}
    }
    \caption{Comparing QR diagrams for (a) very high-order KG solution with Roe-type numerical flux and (b) Energy-preserving KG solution without any additional interface stabilisation in the numerical flux, see Eq.~\eqref{eq:centralKG}, which is stable for $m \leq 4$.}
    \label{fig:X5}
\end{figure}

To confirm that the very high-order limit is representative of an energy conserving dynamics, we compare in Fig.~\ref{fig:X5} the QR profiles of case $m = 16$, $n_{el} = 7$ with that of the kinetic energy preserving KG split form without interface dissipation, cf.\ \eqref{eq:centralKG}, which is stable for $m = 4$. It is possible that such a stable solution is only available because of non-negliable dispersion errors that might be staggering the localisation of the thin shear layers traditionally regarded as connected to inviscid TGV crashes \cite{brachet1992numerical}. As previously noted, the stability of simulations that use split forms without interface dissipation is delicate and crash for $m > 4$, see \cite{gassner2016split}. In any case, as one can see from Fig.~\ref{fig:X5}, the very high order Roe-based QR profile is symmetric and very similar to that of the solution computed with a ``fully central'' spilt form DG scheme without additional dissipation at interfaces. Clearly, both solutions are highly affected by over-energetic bump-related scales which favour equipartition of energy and a symmetrical distribution of kinematic states. The QR profile of the LLF-based solution for case $m = 16$, $n_{el} = 7$ (not shown) is similar to those of Fig.~\ref{fig:X5}.

\begin{figure}[!ht]
	\centering
	\subfloat[$m=3$, $n_{el}=37$: Roe numerical flux.]
	{
		\includegraphics[scale=0.60]{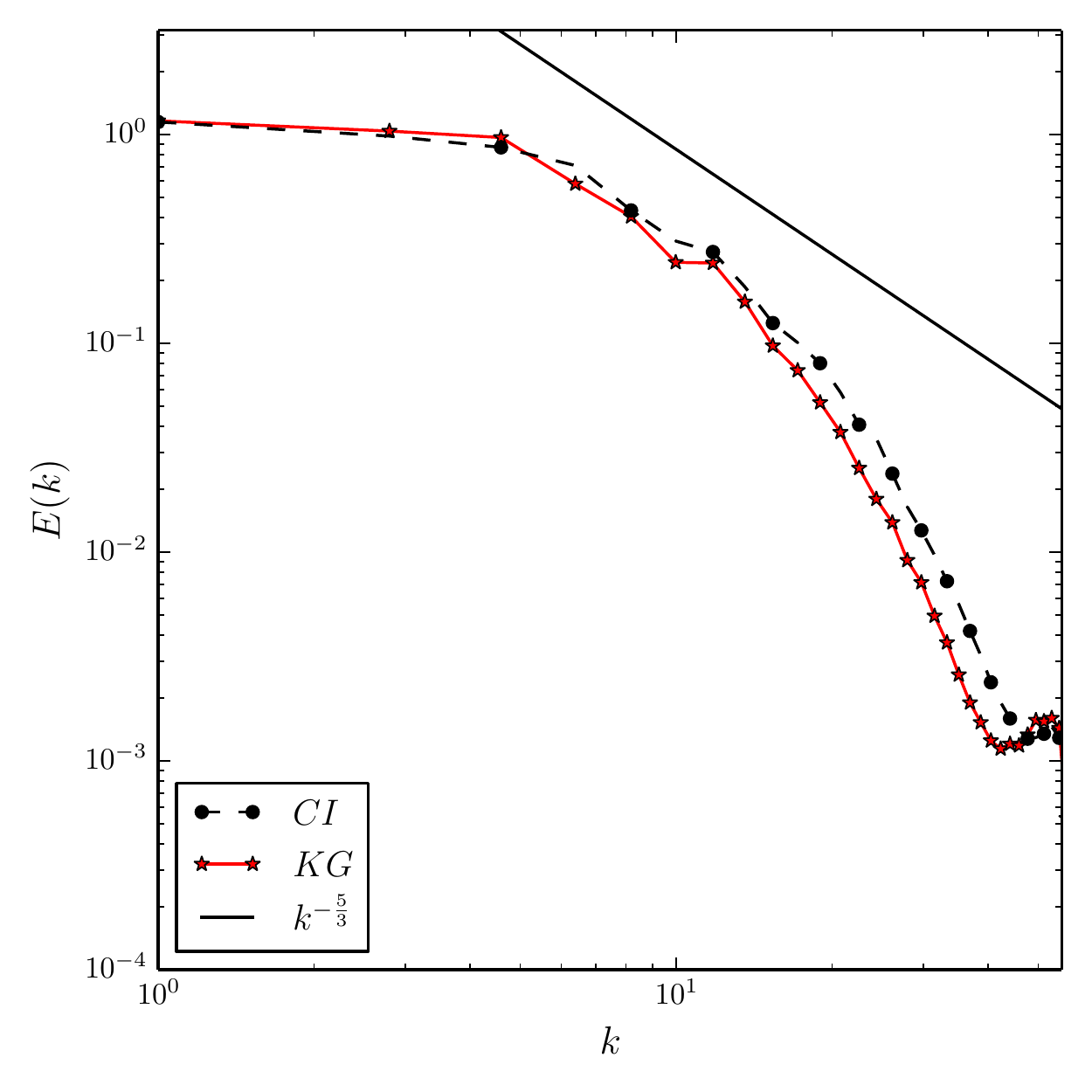}
	}
	\subfloat[$m=3$, $n_{el}=37$: LLF numerical flux.]
	{
		\includegraphics[scale=0.60]{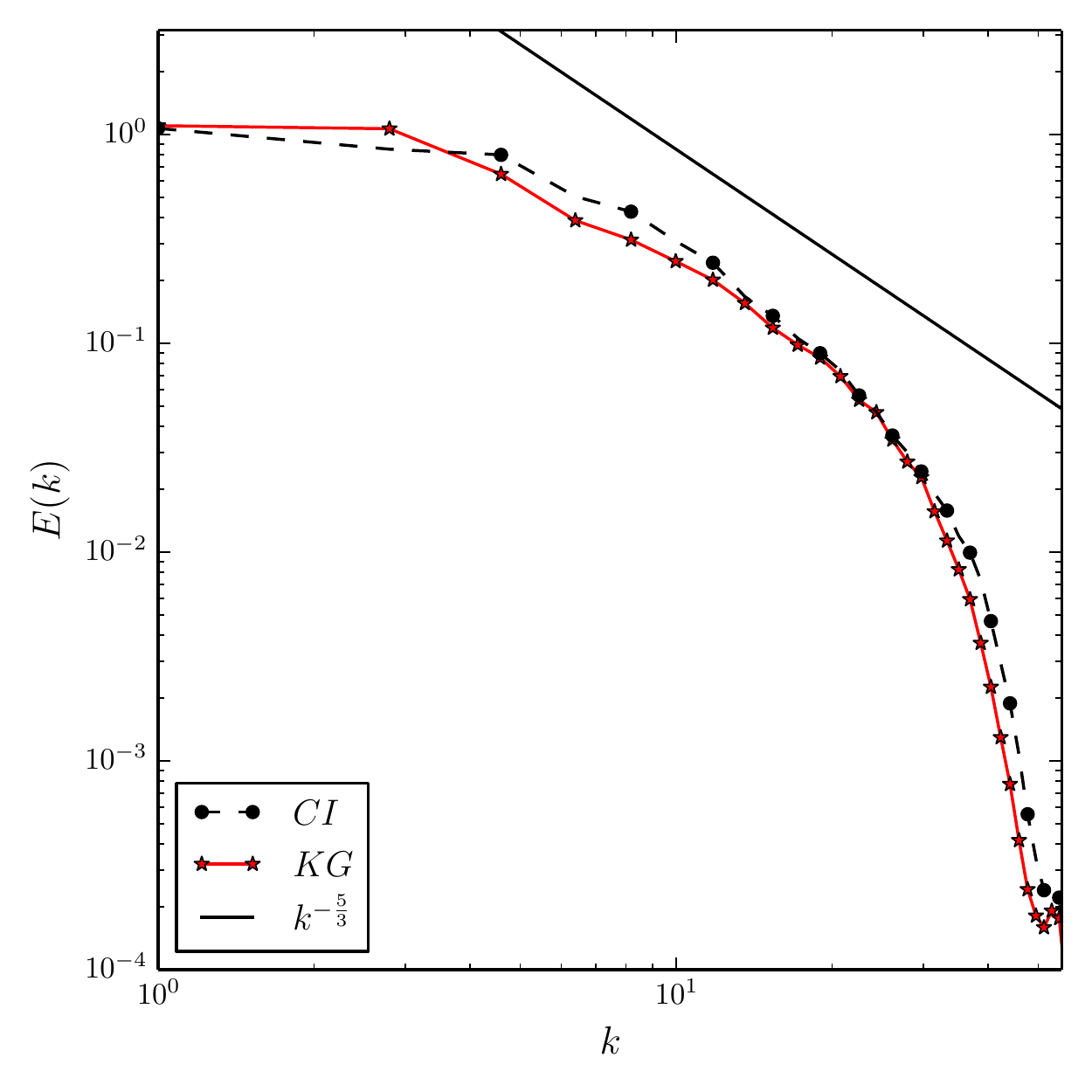}
	}
	\caption{Comparison between energy spectra of consistent integration and split forms for case $m = 3$, $n_{el} = 37$.}	\label{fig:X6}
\end{figure}

Finally we consider some results obtained with low-order discretisations, namely, $m \leq 3$. In Fig.~\ref{fig:X6}, we show the $m = 3$ energy spectra obtained both for split forms and consistent integration. As one can see, non-negligible differences are present. We stress that LLF results seem bump-free, but actually the additional small-scale energy of LLF spectra (when compared to Roe spectra) might be related to a mild energy bump. This is supported by the evolution of the Roe and LLF energy spectra at later times (not shown), when a small-scale spectral region seems to persist with notably high energy in the LLF-based spectrum during the third phase of the TGV flow (nearly homogeneous decay).

\begin{figure}[!ht]
	\centering
	\subfloat[$m = 2$, $n_{el} = 56$: CI LLF numerical flux.]
	{
		\includegraphics[scale=0.505]{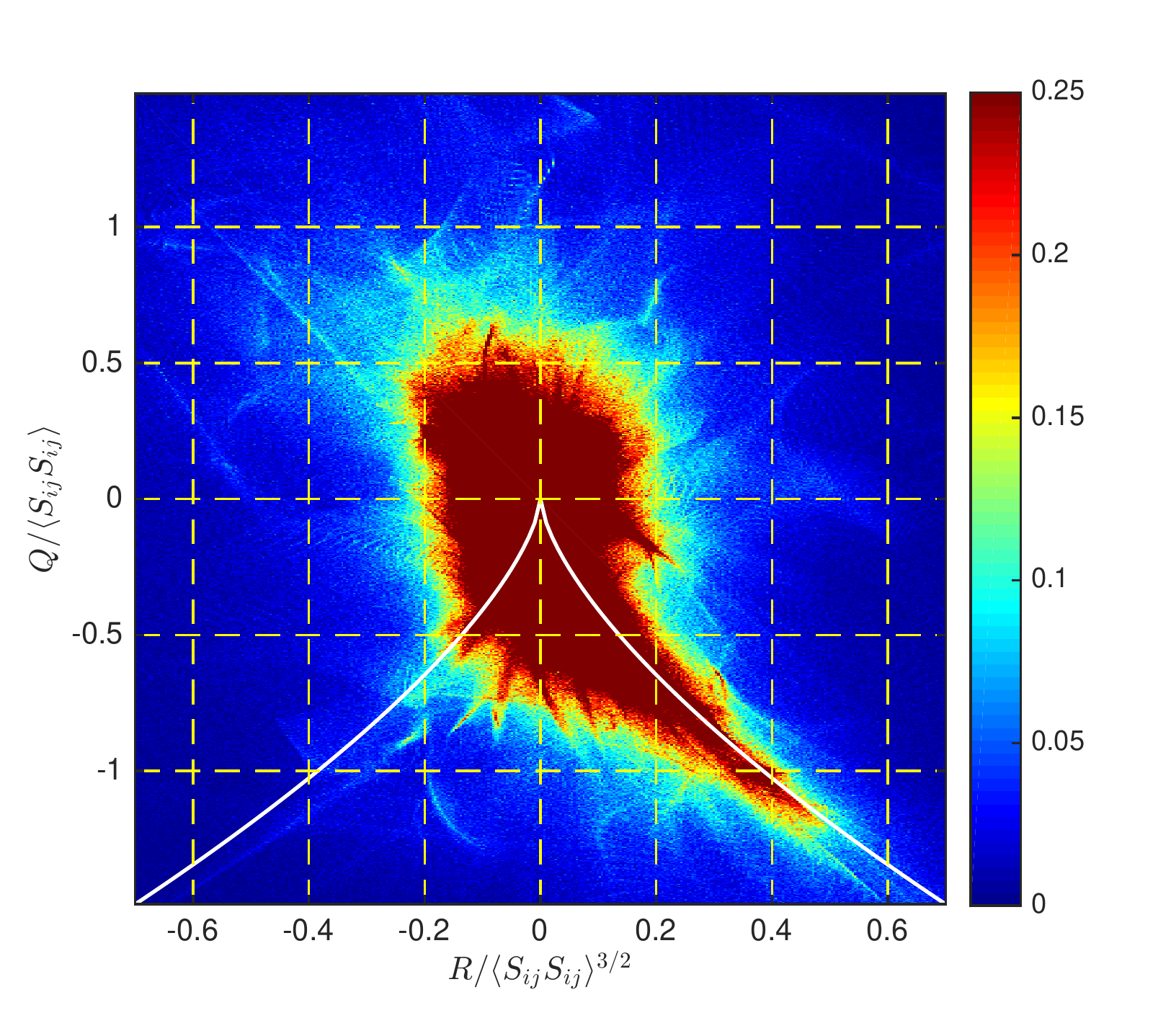}
	}
	\subfloat[$m = 2$, $n_{el} = 56$: CI Roe numerical flux.]
	{
		\includegraphics[scale=0.505]{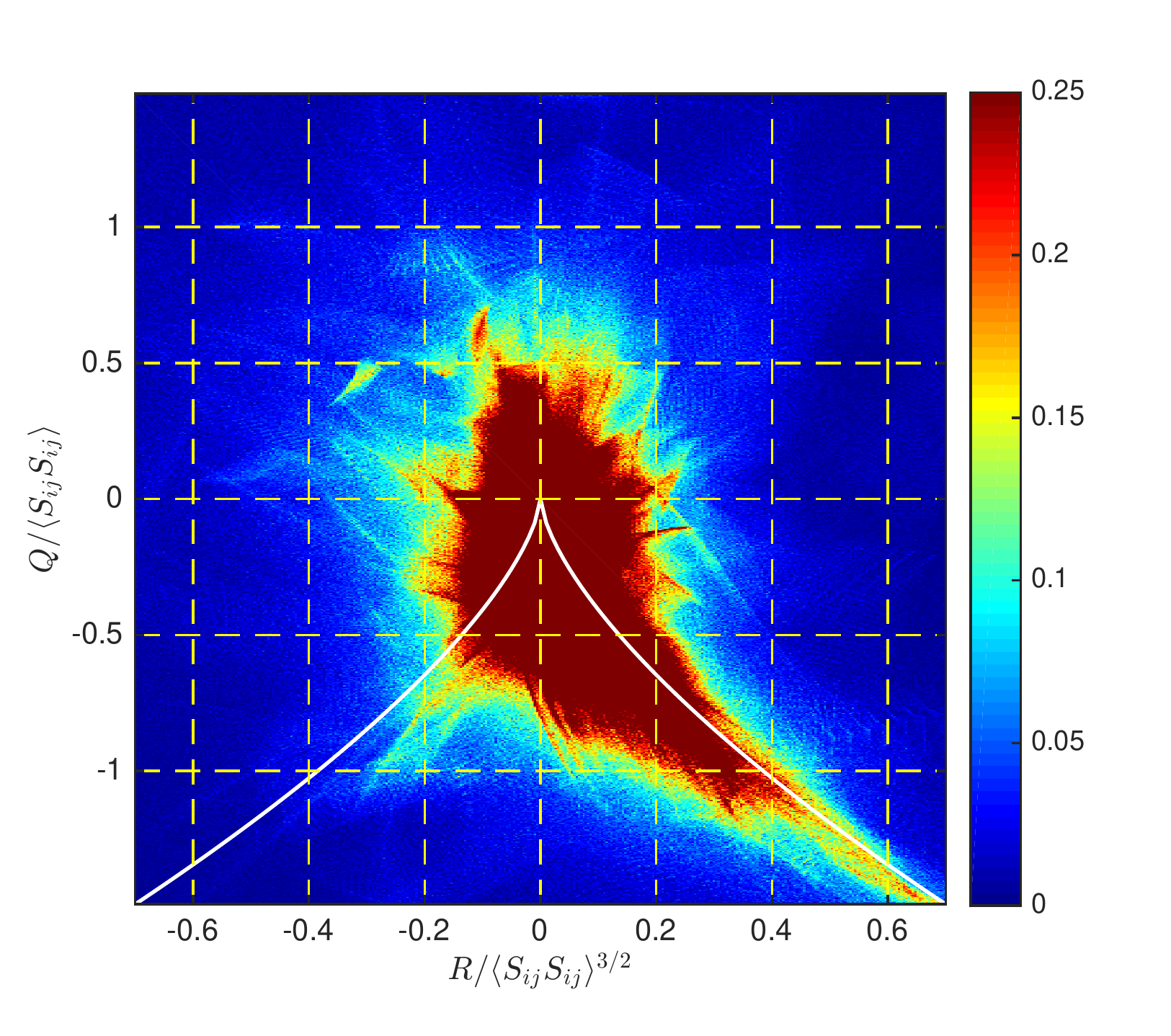}
	}
	\\
	\subfloat[$m = 2$, $n_{el} = 56$: KG LLF numerical flux.]
	{
		\includegraphics[scale=0.505]{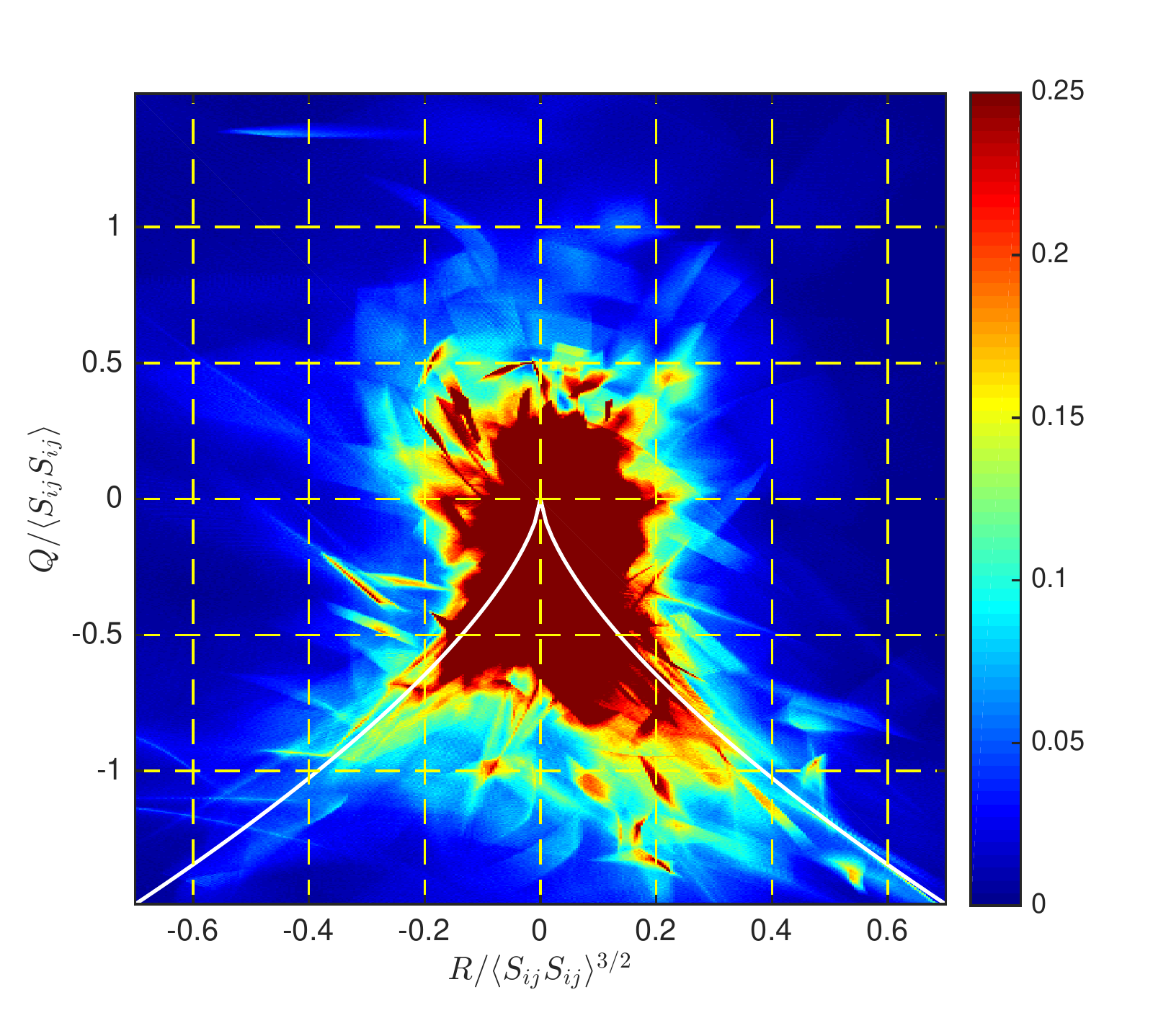}
	}
	\subfloat[$m = 2$, $n_{el} = 56$: KG Roe numerical flux.]
	{
		\includegraphics[scale=0.505]{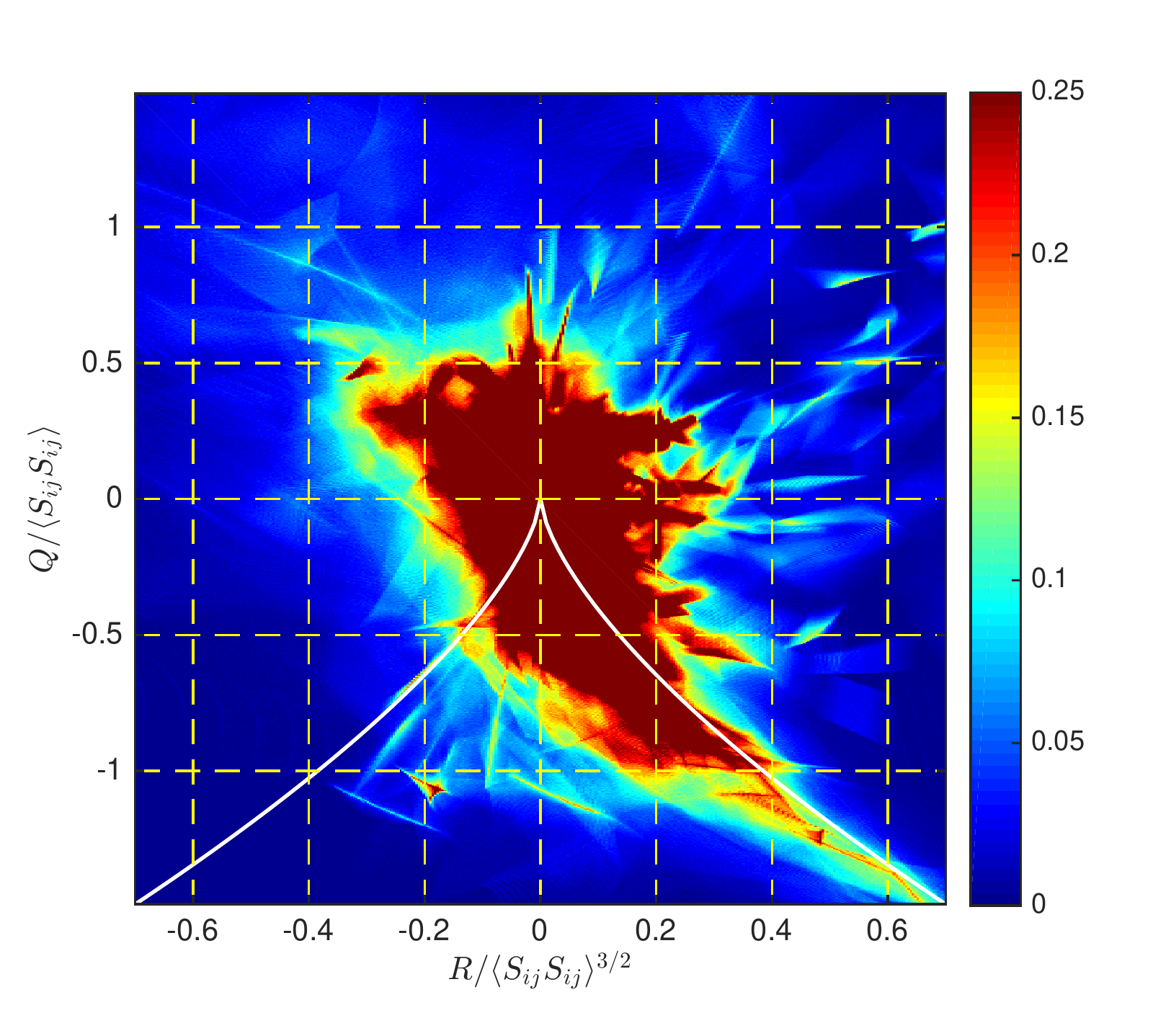}
	}
	\caption{QR diagrams of case $m = 2$, $n_{el} = 56$ at $t = 9$ for consistent integration (top) and split forms (bottom) with LLF (left) and Roe-type (right) numerical fluxes.}
	\label{fig:X7}
\end{figure}

Additional insight on the quality of low-order cases can be gained from QR diagrams, which are given in Fig.~\ref{fig:X7}, now, for $m = 2$. These confirm that, at sufficiently low orders, split forms and consistent integration results can differ. More importantly, they show that the turbulent kinematics of low-order DGSEM solutions are not as clean and accurate when compared to solution at (moderately) higher orders for the same DOFs. However, the overall shape of the low-order QR profiles is still reasonably correct, and these seem even more physical than those obtained by discretisations of very high order.

\section{Conclusions}
\label{sec:conclusions}

In this work we described two strategies to dealias the computation of under-resolved turbulent flows using the DGSEM framework. The first was the well-known technique of consistent/over- integration (sometimes referred to as polynomial dealiasing) where the approximation of variational forms in the DG framework are enriched with additional quadrature points to remove aliasing errors associated with non-linear terms. The other dealiasing strategy was to reformulate the PDE into an equivalent, at the continuous level, split form expression that employs an average of the conservative and non-conservative forms of the equation. Relevant DGSEM split forms can however be shown to remain conservative. It is important to note that both consistent/over- integration and the split form DG schemes require additional computational effort compared to the standard DGSEM. The exploration of optimising implementations of each dealiasing strategy is a subject of our ongoing research.

A preliminary investigation about the built-in dealiasing mechanism of DGSEM split forms has been conducted in one dimension through a frozen Burgers' turbulence scenario. This analysis showed that, without consistent/over-integration, standard conservative DG formulations tend to overestimate the energy content of high-order polynomial modes, whereas advective formulations tend to under-estimate them. Due to the averaging of these two opposite tendencies, the relevant quadratic split form DGSEM was shown to balance aliasing errors and to keep energy content slightly below the correct levels. This latter feature is believed to also help improve robustness in under-resolved computations by keeping the growth of energy levels ``under control''.

Subsequently, the robustness of DGSEM approaches based on consistent integration and split forms for the inviscid TGV problem was discussed in detail. This brought together results from two recent works \cite{gassner2016split,moura2016ontheeddy} \comment{as well as new solutions obtained with lower polynomial orders} and clarified how the observed TGV instabilities can be either induced or suppressed by the underlying numerics. Although TGV instabilities are not entirely understood at this point, it was argued that discretisations that promote an energy-conserving bias for the TGV flow are prone to develop instabilities. In any case, relevant split forms clearly display superior robustness and only lack stability when central fluxes are used in place of a Riemann solver. \comment{Still, the novel part of the TGV analysis concerned the quality of the turbulent solutions obtained with stable split forms. This assessment was conducted by comparing these to solutions obtained with consistent/over-integration, which, as the standard dealiasing technique, is believed to yield the best accuracy DG-based model-free turbulence computations can offer.}

In the accuracy assessment of consistent integration and split forms for the inviscid TGV flow, kinetic energy spectra and QR diagrams have been analysed. At moderately high orders (i.e.\ around sixth order), the performance of relevant split forms was shown to be very close to that obtained with consistent integration. In this case, energy spectra and QR diagrams both exhibited the expected physical trends. At lower orders (i.e.\ third order or below), however, the two DGSEM approaches exhibited non-negligible differences and inferior solution quality (for the same number of DOFs). Most surprisingly, at very high orders (i.e.\ eighth order or above), stable solutions showed unphysical features attributed to the energy-conserving bias induced by a sharper dissipative behaviour in wavenumber space. The latter has been found to be especially prominent for the local Lax-Friedrichs flux, owing to its over-upwind character at low Mach numbers. \comment{Although this energy-conserving bias had already been conjectured in a previous short note \cite{moura2016anlesICOSAHOM}, new and stronger evidence for it was presented here.}

In summary, DG-based model-free turbulence computations at very high Reynolds numbers are more likely to yield stable and accurate solutions when complete Riemann solvers (such as Roe's flux) are employed at moderately high orders. For challenging simulations, the use of dealiasing strategies is of key importance. This work highlighted the role of high-order split form discretisations as a robust alternative compared to consistent/over- integration. More importantly, it showed that high-order solutions from both strategies were very similar in terms of solution quality. 
An interesting theme left for subsequent studies concerns the performance of split form DGSEM for cases involving curved elements, whose type of aliasing error, i.e. geometric aliasing, has a different source than that of under-resolved turbulence.

\section*{Acknowledgement}
Andrew Winters was partially supported by the ``Mobility Grant for National and International Young Faculty'' from the University of Cologne. Stefanie Walch thanks the Deutsche Forschungsgemeinschaft (DFG) for funding through the SPP 1573 ``The physics of the interstellar medium'' and the funding from the European Research Council via the ERC Starting Grant  ``The radiative interstellar medium'' (\texttt{RADFEEDBACK}). Gregor Gassner has been supported by the European Research Council (ERC) under the European Union's Eights Framework Program Horizon 2020 with the research project \textit{Extreme}, ERC grant agreement no. 714487.

Rodrigo Moura would like to acknowledge funding under the Brazilian Science without Borders scheme. Spencer Sherwin additionally acknowledges support as Royal Academy of Engineering Research Chair under grant 10145/86. 

The simulations in this work were partially performed using the Cologne High Efficiency Operating Platform for Sciences (CHEOPS) HPC cluster at the Regionales Rechenzentrum K\"{o}ln (RRZK), University of Cologne, Germany. Use of the HPC facilities at Imperial College London is also acknowledged.

\section*{References}

\bibliography{References}

\appendix

\section{Comparison between representative split forms}
\label{sec:compSplit}

Here we explicitly compare two split forms, namely DU and KG. Among various possible split forms, these two are considered well-balanced choices when different computational aspects are taken into account \cite{gassner2016split}. The results discussed in the following were obtained from case $m = 5$, $n_{el} = 23$ with Roe's numerical flux, but other test cases considered (not shown) clearly support the conclusions below.

\begin{figure}[!ht]
	\centering
	\subfloat[Evolution of the enstrophy.]
	{
		\includegraphics[scale=0.60]{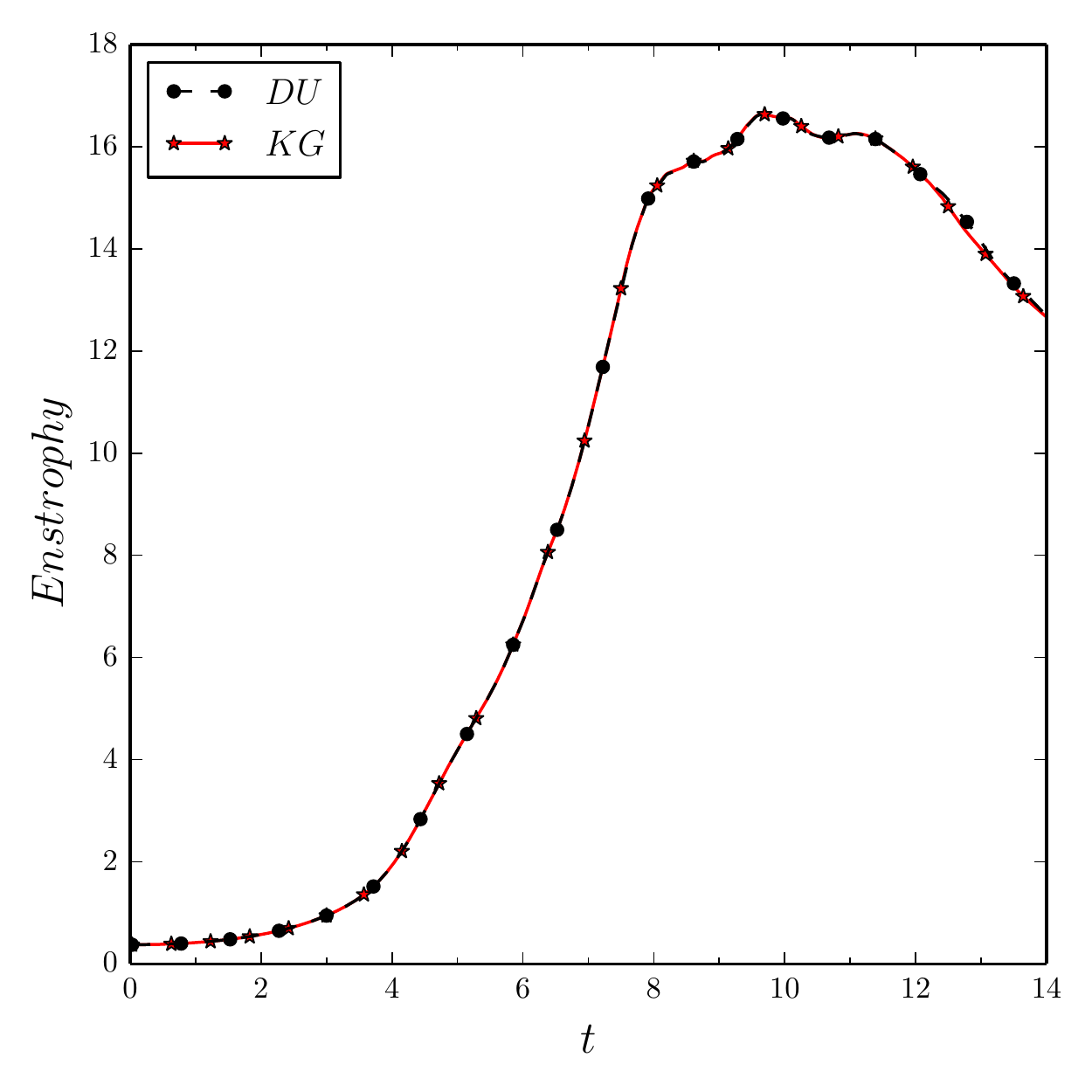}
	}
	\subfloat[Energy spectrum at $t=9$.]
	{
		\includegraphics[scale=0.60]{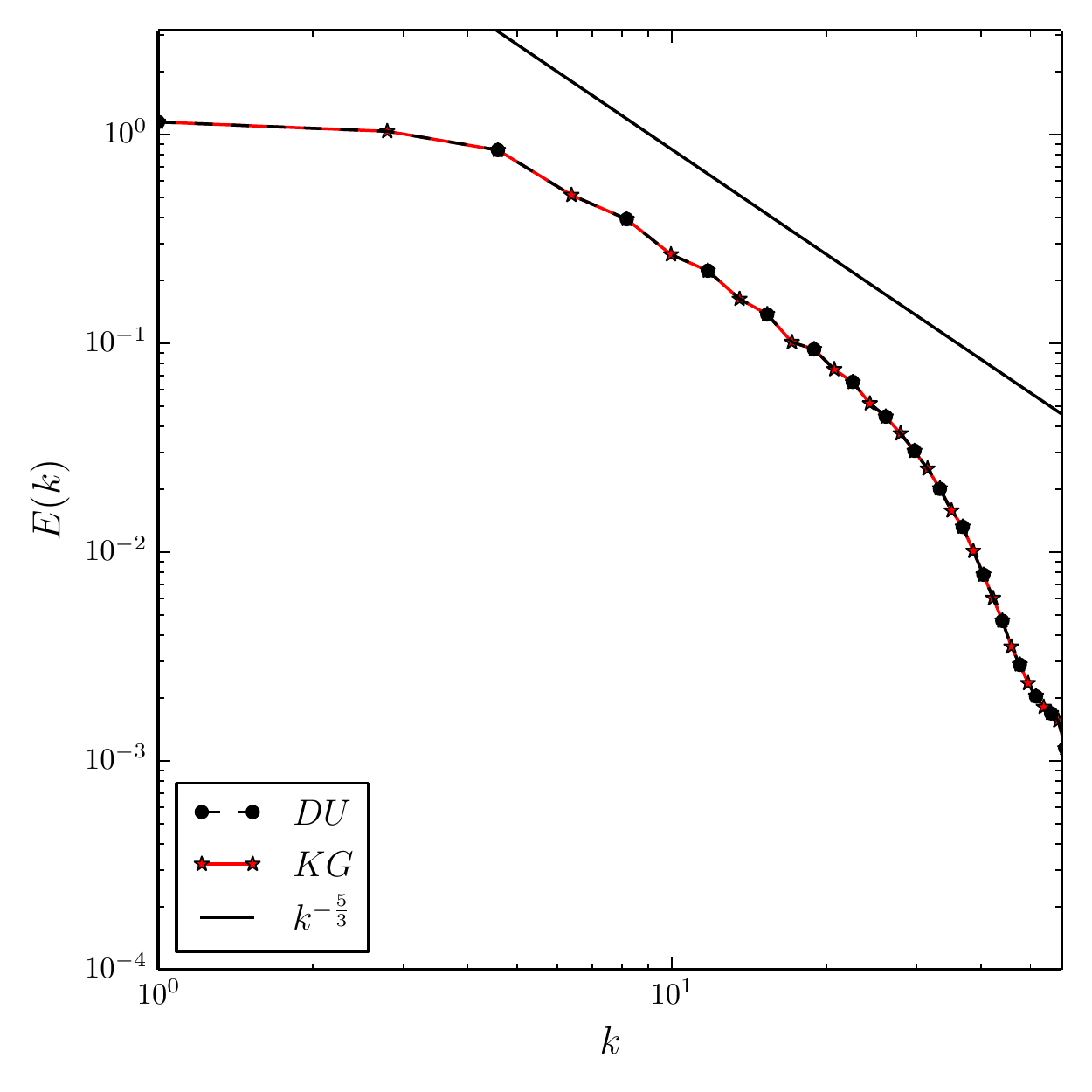}
	}
	\\
	\subfloat[QR diagram for KG split form at $t=9$.]
	{
		\includegraphics[scale=0.505]{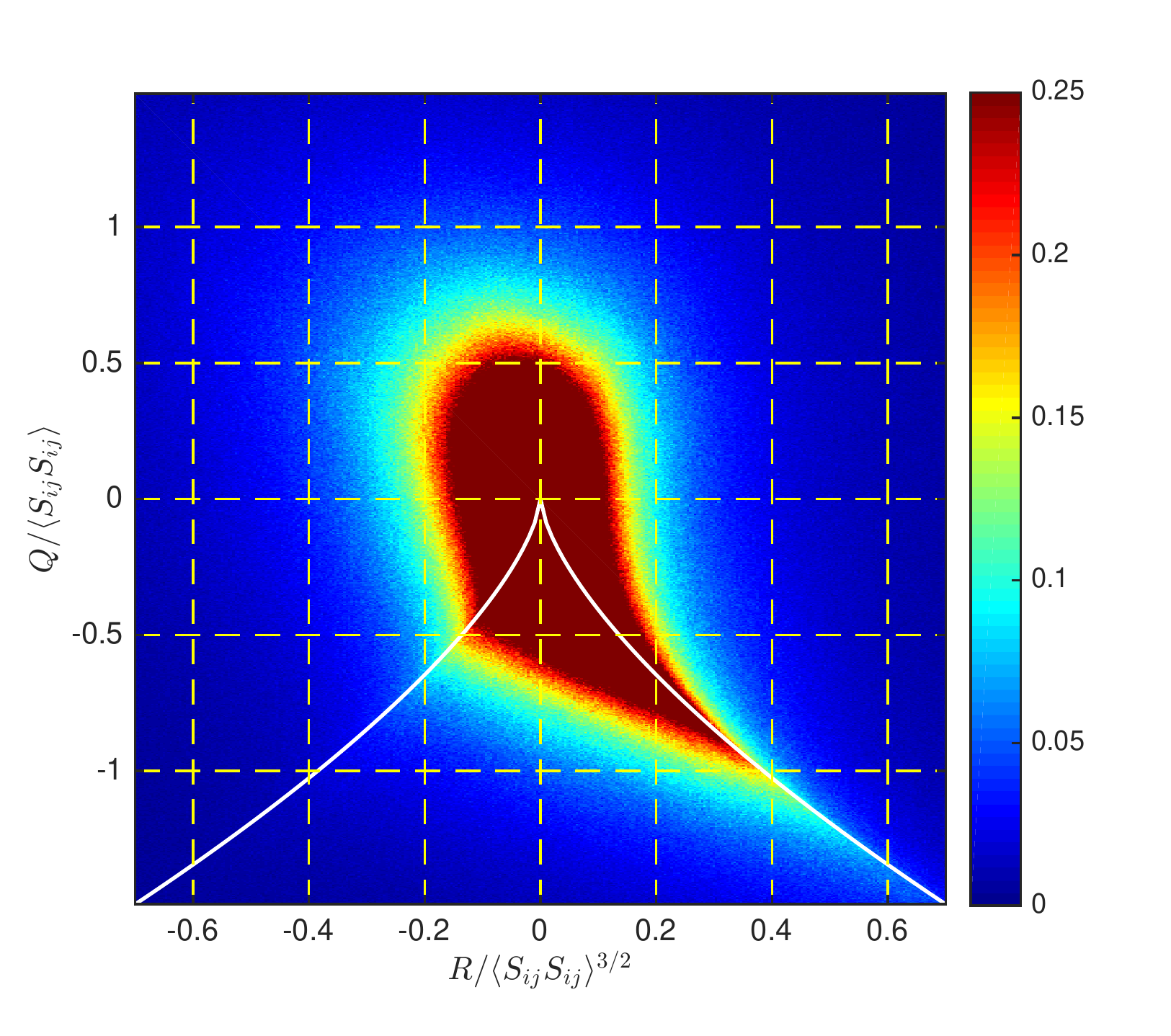}
	}
	\subfloat[Pointwise difference between KG and DU.]
	{
		\includegraphics[scale=0.505]{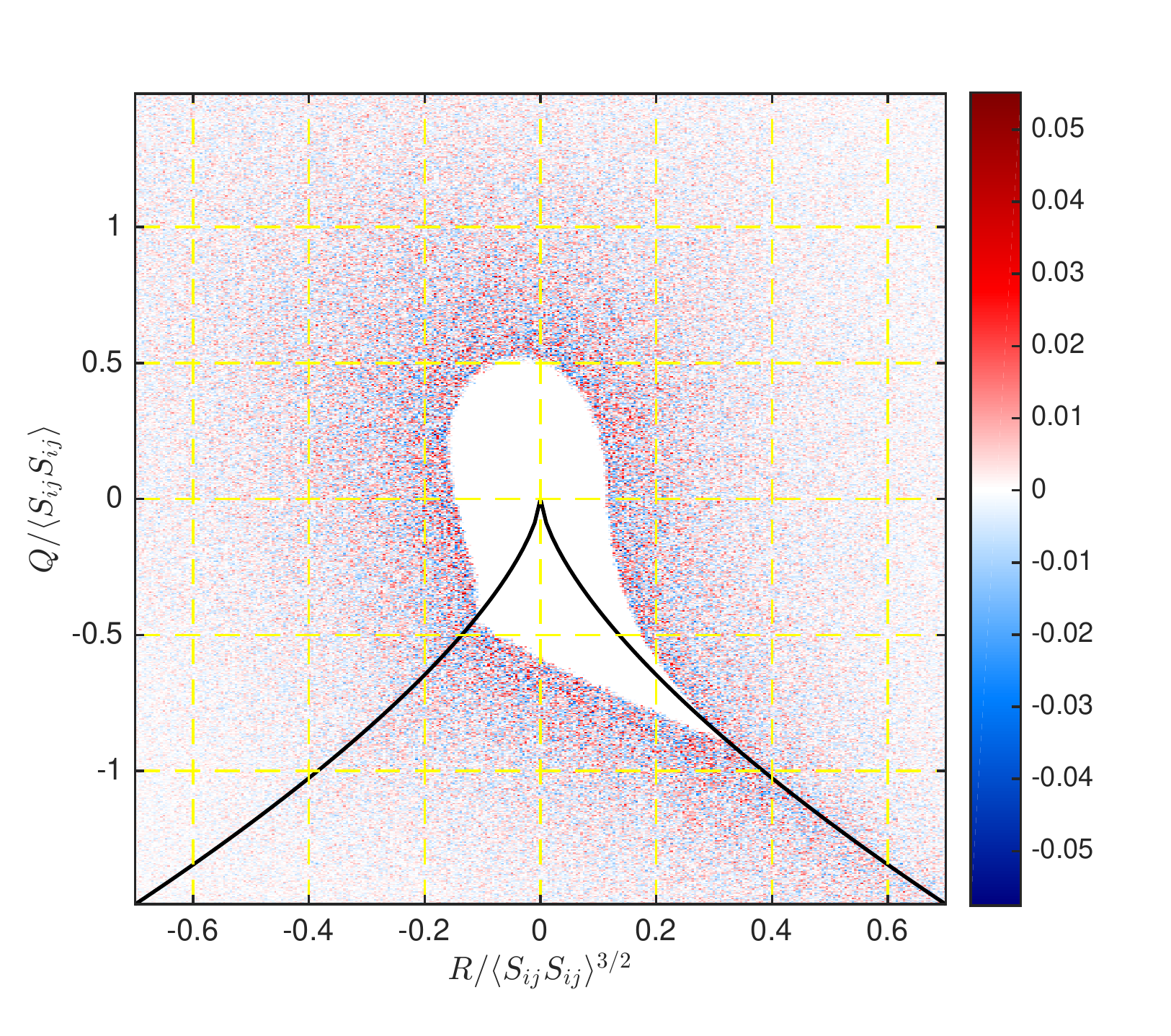}
	}
	\caption{Comparison of split forms DU and KG for the test case $m = 5$, $n_{el} = 23$ with Roe-type numerical flux.}	\label{fig:X1}
\end{figure}

We first consider the evolution of enstrophy, shown in Fig.~\ref{fig:X1}(a), and defined for the TGV box $\Omega$ as
\begin{equation}
\zeta = \frac{1}{|\Omega|} \int_{\Omega}{\, \frac{\, \rho \, |\omega|^2 \!}{2} \, d\Omega} \mbox{,}
\end{equation}
where $|\omega|^2 = \vec{\omega} \cdot \vec{\omega}$ and $\vec{\omega}$ is the vorticity vector. The results for DU and KG are practically indistinguishable. This is also shown to be true in Fig.~\ref{fig:X1}(b) for the energy spectrum at $t = 9$, which is nearly when the peak kinetic energy dissipation is attained and where one expects Kolmogorov's -5/3 slope to be followed \cite{brachet1991direct}. The curves shown correspond to standard three-dimensional energy spectra and have been generated as described in \cite{moura2016ontheeddy}: by probing the solution on a Cartesian set of equispaced points centred in a way that avoids probing data at elemental interfaces. The total number of points used matches $N_{dof} = (n_{el} \, m)^3$ for each test case.

The QR diagram of the flow at $t=9$ is given in Fig.~\ref{fig:X1}(c) for the KG split form. We see that this split form DGSEM recovers the teardrop shape that is expected for turbulent flows. Additionally, the KG and DU diagrams were nearly identical and therefore only the KG result is shown here. Nevertheless, we provide a pointwise comparison between the two split form QR diagrams in Fig. \ref{fig:X1}(d).

While these two forms are expected to differ for flows with stronger compressibility effects \cite{gassner2016split}, the results in Fig. \ref{fig:X1} encouraged us to consider only one of the two forms in the analyses in Section \ref{sec:solQuality}. The KG form has been chosen therefore as the representative split form in this study, since it is known to be more robust for properly compressible flows \cite{gassner2016split}.

\section{Discussion of QR diagrams}\label{app:QR}

We provide a brief discussion, motivation and construction of the QR diagram. The QR diagram provides an interesting statistical representation for possible linear local flow trajectories of turbulent kinematics, complete details are given by Chong et al. \cite{Chong1990}. The flow trajectories for incompressible flows are categorized in the space of the joint PDFs of the second (Q) and third (R) invariants of the velocity gradient tensor 
\begin{equation}\label{eq:tensor}
A_{ij} = \frac{\partial u_i}{\partial x_j},
\end{equation}
for a given flow field. The velocity gradient tensor \eqref{eq:tensor} can be decomposed into a symmetric part
\begin{equation}\label{eq:symmetric}
S_{ij} = \frac{1}{2}\left( \frac{\partial u_i}{\partial x_j} +  \frac{\partial u_j}{\partial x_i}\right),
\end{equation} 
and an anti-symmetric part
\begin{equation}\label{eq:antiSymmetric}
W_{ij} = \frac{1}{2}\left( \frac{\partial u_i}{\partial x_j} -  \frac{\partial u_j}{\partial x_i}\right).
\end{equation}
The symmetric part \eqref{eq:symmetric} is defined as the strain rate tensor and the anti-symmetric part \eqref{eq:antiSymmetric} as the rotation rate tensor. The second and third invariants of \eqref{eq:tensor} are given by
\begin{equation}\label{eq:invariants}
\text{Q} = -\frac{1}{2}A_{ij}A_{ji},\quad \text{R} = -\frac{1}{3}A_{ij}A_{jk}A_{ki},
\end{equation}
respectively. Note in the QR diagrams, e.g. Fig.~\ref{fig:X1}(c), we normalise Q with $\langle S_{ij}S_{ij}\rangle$ and R by $\langle S_{ij}S_{ij}\rangle^{3/2}$, where $\langle\cdot\rangle$ is the $L^2$ inner product. 
The white curve in the QR diagrams corresponds to the zero discriminant $\frac{27}{4}\text{R}^2+\text{Q}^3=0$ that defines a boundary between topologically distinct flow patterns \cite{Chong1990}. Thus, one can interpret crossing the discriminant curve as a change in the flow trajectory of turbulent kinematics. \comment{The canonical teardrop shape observed in QR diagrams extracted from the DNS of isotropic turbulence, see e.g.\ \cite{Laizet2015}, very much resembles that shown in Fig. A.8(c), where the characteristic "tail" aligns with the zero discriminant curve for $R>0$.}

\end{document}